\documentclass[leqno,12pt]{article}%
\usepackage{amsmath,amssymb}
\usepackage{xparse}
\usepackage{amsthm}
\usepackage{dsfont}
\usepackage{amsfonts}%
\usepackage{color}

\newtheorem{Theorem}{\hspace{\parindent}\bf Theorem}[section]
\newtheorem{Lemma}{\hspace{\parindent}\bf Lemma}[section]

\newtheorem{Corollary}{\hspace{\parindent}\bf Corollary}[section]

\begin{document}

\title{\textbf{Propagation of solutions of the Porous Medium Equation with reaction and their travelling wave behaviour}}
\author{by \\ Alejandro G\'arriz}

\maketitle

\

\begin{abstract}
We consider reaction-diffusion equations of porous medium type, with different kind of reaction terms, and nonnegative bounded initial data. For all the reaction terms under consideration there are initial data for which the solution converges to 1 uniformly in compact sets for large times. We will characterize for which reaction terms this happens for all nontrivial nonnegative initial data, and for which ones there are also solutions converging uniformly to 0. Problems in this family have a unique (up to translations) travelling wave with a finite front and we will see how its speed gives the asymptotic velocity of all the solutions with compactly supported initial data. We will also prove in the one-dimensional case that solutions with bounded compactly supported initial data converging to 1 do so approaching a translation of this unique traveling wave. We will prove a similar result for non-compactly supported initial data in a certain class.
\end{abstract}



\vskip 4cm

\noindent{\makebox[1in]\hrulefill}\newline2010 \textit{Mathematics Subject
Classification.}
35K55. 
35B40. 
35K65, 
76S05 
\newline\textit{Keywords and phrases.} Reaction-diffusion, porous medium equation, travelling wave behaviour.

\newpage

\section{Introduction and main results}
\label{sect-introduction} \setcounter{equation}{0}

The aim of this paper it to characterize the large time behaviour of solutions to the Cauchy problem
\begin{equation}\label{eq:main}
u_t = \Delta u^m + h(u)\quad\text{in }Q:=\mathbb{R}^N\times\mathbb{R}_+,  \qquad u(\cdot,0)=u_0 \geq 0 \quad\text{in }\mathbb{R}^N,
\end{equation}
where $m>1$. The nonlinearity $h$ is assumed to be in $C^1(\overline{\mathbb{R}}_+)$ and to fulfill, for some $a\in[0,1)$,
\begin{equation}\label{eq:reaction}
h(0)=0, \quad
h(u)\leq 0\text{ if } u \in [0,a],\quad h(u)>0\text{ if } u \in (a,1),\quad h(u)<0\text{ if } u > 1.
\end{equation}

Let us also ask our reaction term to be \lq\lq exponentially stable" in $u=1$, that is,
\begin{equation}\label{eq:reaction_deriv_in_1}
h^\prime (1)  <0.
\end{equation}

We are mainly interested in the case of bounded and compactly supported initial functions $u_0$, though in the one-dimensional scenario we will deal also with a different class of initial data.

If $a=0$, we are in the so-called monostable case, which includes the Pearl-Verhulst (or logistic) reaction nonlinearity $h(u)=u(1-u)$. The case $a>0$ contains as particular instances the classical bistable nonlinearities, when $h(u)<0$,  and the combustion ones, when $h(u)=0$ for $u\le a$ (in this case $a$ is the ignition temperature).

Problem~\eqref{eq:main}--\eqref{eq:reaction} may be used, for example, to describe the growth and propagation of a spatially distributed biological population whose tendency to disperse depends on the population density; see for instance~\cite{Gurney-Nisbet-1975, Gurtin-MacCamy-1977}. It can also be regarded as a generalization of the semilinear case $m=1$. When $a=0$, this semilinear problem  was introduced to study spreading questions in population genetics independently by Fisher~\cite{Fisher-1937},  and by Kolmogorov, Petrovsky and Piscounov~\cite{Kolmogorov-Petrovsky-Piscounov-1937}; see also the work of Sanchez-Gardu\~no and Maini~\cite{SanchezGarduno-Maini-1994} for more applications. When $a>0$, it was proposed by Zel'dovich to study combustion problems~\cite{Zeldovich-1948}; see also~\cite{Kanel-1964}. Finally, in the field of astronomy, our equation was proposed by Newman and Sagan as a model for the propagation of intergalactic civilizations~\cite{Newman-Sagan-1981}.

The equation in~\eqref{eq:main} is degenerate when $u=0$, and does not have in general classical solutions.  A function $u$  is a \emph{weak solution} to problem~\eqref{eq:main} if $u,\nabla u^m, h(u)\in L^1_{\rm loc}(Q)$, and
\begin{equation}
\label{eq:weak.solution}
\int_{\mathbb{R}^N}u_0\varphi(0)+\int_0^t\int_{\mathbb{R}^N}
\big(u\varphi_t-\nabla u^m\cdot\nabla\varphi+h(u)\varphi\big)=0
\end{equation}
for each $\varphi\in C_{\rm c}^\infty(\overline Q)$. If the equality in~\eqref{eq:weak.solution} is replaced by \lq\lq$\le$'' (respectively, by \lq\lq$\ge$'') for $\varphi\ge 0$, we have a subsolution (respectively, a supersolution). Problem~\eqref{eq:main} has a unique bounded weak solution when the initial function is bounded. This solution is continuous in $Q$, and satisfies in addition that $\nabla u^m\in L^2_{\rm loc}(Q)$. There is moreover a comparison principle between bounded subsolutions and supersolutions~\cite{Sacks-1983,dePablo-Vazquez-1991}.

The large time behaviour will be given in terms of \emph{travelling wave} solutions (abreviated TW in what follows). By this we mean a solution of the form $\bar u_c(x,t)=V_c(x-ct)$ for some \emph{speed} $c$ and  \emph{profile}  $V_c$ (depending on $c$), which should satisfy $V_c\in C(\mathbb{R})$, $(V^m_c)^\prime,h(V)\in L^1_{\rm loc}(\mathbb{R})$,  and
$$
\displaystyle\int_{\mathbb{R}} \left((V_c^m)^\prime\varphi^{\prime} +c\ V_c \ \varphi^\prime - h(V_c)\varphi\right) =0
$$
for all $\varphi \in C_{\rm c}^{\infty} (\mathbb{R})$.

Two wave profiles are said to be indistinct if one is a translation of the other one. A monotonic travelling wave solution is called a \emph{change of phase type} or a \emph{wavefront from 1 to 0} if it connects the two equilibrium states 1 and 0, that is,
$$
\lim\limits_{\xi \rightarrow - \infty} V_{c}(\xi) = 1,\qquad\lim\limits_{\xi \rightarrow  \infty} V_{c}(\xi) = 0.
$$

If the reaction~\eqref{eq:reaction} falls in the case $a=0$,  there exists a minimal speed $c^* = c^*(m,h)>0$ such that for all $c \geq c^*$ equation~\eqref{eq:main}  has an unique distinct monotonic change of phase type TW profile $V_c$.  There are no admissible TW solutions for $c < c^*$ though.

On the other hand if $a>0$ there is a unique speed $c^*$ for which the equation has a wavefront. In this case the sign of $c^*$ matches the sign of
$\displaystyle\int_0^1 h(u)u^{m-1}\ du$. The asymptotic behaviour when $c^*\le 0$ has already been considered in~\cite{Hosono-1986}. Hence, in the sequel we will always assume
\begin{equation}
\label{eq:cond.Hosono}
\int_0^1 h(u)u^{m-1}\ du>0.
\end{equation}

\normalcolor

The TW with speed $c^*$ satisfies $V_{c^*}<1$ and is finite, i.e., there exists a value $\xi_0$ such that $V_{c^*}(\xi) =0$ for all $\xi \geq \xi_0$ and $V_{c^*}(\xi) >0$ for all $\xi < \xi_0$. Moreover, $V^\prime_{c^*}<0$ for all $\xi < \xi_0, V_{c^*}^m \in C^1(\mathbb{R})$ and $(V_{c^*}^m)^\prime (\xi) \rightarrow 0$ as $\xi \rightarrow \pm\infty$.


Finally, $V_{c^*}$ satisfies that
\begin{equation}\label{profile:property}
\lim\limits_{\xi \rightarrow \xi_0^-} m\ V_{c^*}^{m-2}(\xi)\ V_{c^*}^\prime(\xi) = \lim\limits_{\xi \rightarrow \xi_0^-} \left(\frac{m}{m-1}\ V_{c^*}^{m-1}(\xi)\right)^{\prime} = -c^*.
\end{equation}
This is a hint of the importance of the \emph{pressure} variable $p=\frac{m}{m-1}u^{m-1}$.

From now on, $V_{c^*}$ will denote the TW profile with $\xi_0=0$. Further information about these results and many others can be found in~\cite{Gilding-Kersner-2004} and the references therein.



 In the first part of this work we will generalize most of the results for $m=1$ presented by Aronson and Weinberger in~\cite{Aronson-Weinberger-1978} to the case of density dependent diffusion, i.e., $m>1$. Given an initial data we will study whether the solution converges to 1, a situation that, following the literature, we will denote as \emph{spreading},  or to 0, which we will call \emph{vanishing}. Our first result shows that for all the reaction nonlinearities that we are considering there are initial data for which spreading happens.

\begin{Theorem}
There exists a three parameter family of continuous, bounded and compactly supported functions $v(x; x_0,\eta,\rho)$ (see Section~\ref{sect-propagation.vs.vanish} for a precise description) such that  if
$$
u(x,0) \geq v(x; x_0,\eta,\rho)
$$
for some $x_0 \in \mathbb{R}^N, \eta, \rho >0$, then $u$ converges to 1 uniformly on compact sets.
\end{Theorem}

It turns out that for certain monostable nonlinearities $h$ nontrivial solutions always spread, independently of the mass of the initial datum. Following~\cite{Aronson-Weinberger-1978}, this is named as the \textit{hair-trigger effect}.

\begin{Theorem}
\label{thm:spreading.vs.vanishing}
Let $h$ satisfy~\eqref{eq:reaction} with $a=0$ and
$$
\displaystyle\liminf\limits_{u \rightarrow 0} \frac{h(u)}{u^{m + 2/N}} > 0.
$$
If $u\not\equiv 0$, then $u$ converges to 1 uniformly on compact sets.
\end{Theorem}

Let us remark that if $h(u)\leq ku^p$, with $p\geq m+ 2/N$ and $k$ a positive constant, comparison with the problem with reaction $ku^p$ shows that for certain small initial data our solution asymptotically vanishes; see for instance~\cite{Samarskii-Galaktionov-Kurdyumov-Mikhailov-Book}.


Now we ask ourselves how fast is the spreading. If we move from a point $y_0\in \mathbb{R}^N$ in a certain direction with a slow speed $c$ in the limit we will see only the value $u=1$, while if $c$ is too fast we will surpass the free boundary of our solution, and thus we will only see the value $u=0$. Too slow will translate to $c<c^*$, and too fast to $c>c^*$. In this sense, the speed $c^*$ is called the \emph{critical speed} of the problem.

\begin{Theorem}
If the initial datum is bounded and compactly supported,  given any $y_0\in\mathbb{R}^N$ and $c>c^*$ there is a value $T$ such that
$$
u(y,t) = 0\quad \text{for }|y - y_0| \geq ct,\, t\ge T.
$$
In addition, if spreading happens, for any $c \in (0,c^*)$ we have
$$
\lim\limits_{t \rightarrow \infty} \min\limits_{|y - y_0| \leq ct} u(y,t) = 1.
$$
\end{Theorem}
This question has been recently analyzed in~\cite{Audrito-Vazquez-2017,Audrito-Vazquez-FDE-2017} for a wider class of diffusion operators, though with less general reaction terms. In particular, these papers require $h$ to be concave in the region where it is positive, an assumption that we do not need.

The results of this first part are presented in Section~\ref{sect-propagation.vs.vanish}, after two sections of a more technical nature. In  Section~\ref{sect-preliminaries} we obtain an important estimate for the gradient of the pressure that plays a key role in the study of the asymptotic profile using Bernstein's method instead of the more common nowadays estimates of semiconcavity; see~\cite{Perthame-Quiros-Vazquez-2014,Du-Quiros-Zhou-Preprint}. This is done in this way due to the wide class of reaction functions we are allowing here. In Section~\ref{sect-subsupersolution} we construct two families of sub- and supersolutions that play a ey role in this paper.

The second part of the paper is devoted to study the asymptotic behaviour of solutions when they converge to 1. We manage to give uniform convergence in dimension $N=1$ in moving coordinates  to  a travelling wave profile. Without loss of generality, we may assume that $u_0\in C(\mathbb{R})$.

We start by considering a class of initial data with unbounded support for which the analysis is similar, though slightly simpler, than the one for the class of bounded and  compactly supported initial data.

Thanks to~\eqref{eq:reaction_deriv_in_1} we know that there is a value  $\delta>0$ such that $h^\prime(u)<0$ for $u\in (1-\delta,1+\delta)$. Let $\mathcal{A}$ be the class of nonnegative continuous functions $u_0$ such that $u_0(x) \equiv 0$ for all $x \geq x_0$  for some $x_0 \in \mathbb{R}$, $\liminf_{x \rightarrow -\infty} u_0(x) \in (1-\delta,1+\delta)$. Solutions with initial data in this class have a \emph{right free boundary},
$$
\zeta(t) \equiv \inf \{r \in \mathbb{R} : u(x,t)=0\text{ for all }x \geq r\},
$$
and converge to the travelling wave $V_{c^*}$ both in shape and speed.
\begin{Theorem}
\label{thm:convergence.class.A}
Let $u$ be a weak solution to~\eqref{eq:main} with $u_0 \in \mathcal{A}$, and let $\zeta$ be the function giving its right free boundary. Then there exists a $\xi_0 \in \mathbb{R}$ such that
\begin{equation}\label{eq:convergence.class.A}
\lim\limits_{t \rightarrow \infty} \sup\limits_{x \in \mathbb{R}} |u(x,t) - V_{c^*}(x - c^* t - \xi_0)|= 0,\qquad\lim\limits_{t\to\infty}\zeta(t)-c^*t=\xi_0.
\end{equation}
\end{Theorem}
This theorem will be proved in Section~\ref{sect-tw.class.A}. The case $h(u)=u^p(1-u)$ with $p\in[1,m]$ was already considered in~\cite{Biro-2002}; see also~\cite{Du-Quiros-Zhou-Preprint}.

We remark that one can weaken the hypothesis $\liminf_{x \rightarrow -\infty} u_0(x) \in (1-\delta,1+\delta)$ in certain cases. For example, if the nonlinearity $h$ is such that we have the hair-trigger effect, one only needs that $\liminf_{x \rightarrow -\infty} u_0(x) >0$. Conditions leading to weaker assumptions are presented in Corollary~\ref{coro:lim_inf}.

We arrive to our main result, a description of the large time behaviour for spreading solutions with bounded and compactly supported initial data in dimension $N=1$. Let us describe the left and right free boundary as
$$
\begin{array}{lcc}
\zeta_+(t) \equiv \inf \{r \in \mathbb{R} : u(x,t)=0\text{ for all }x \geq r\}, \\[8pt]
\zeta_-(t) \equiv \sup \{r \in \mathbb{R} : u(x,t)=0\text{ for all }x \leq r\}.
\end{array}
$$

\begin{Theorem}
\label{thm:convergence.compact.support}
Let $u$ be a solution of Problem~\eqref{eq:main} corresponding to a bounded and compactly supported initial data $u_0$ such that $u$ converges to 1 uniformly in compact sets, and let $\zeta_\pm$ be its left and right free boundaries. There exist constants $\xi_\pm\in\mathbb{R}$ such that
\begin{align}
\lim\limits_{t \rightarrow \infty} \sup\limits_{x \in \overline{\mathbb{R}_+}} |u(x+c_*t,t) - V_{c^*}(x-\xi_+)| = 0, \qquad\lim\limits_{t\to\infty}\zeta_+(t)-c^*t=\xi_+,\\
\lim\limits_{t \rightarrow \infty} \sup\limits_{x \in \overline{\mathbb{R}_-}} |u(x-c_*t,t) - V_{c^*}(\xi_--x)| = 0, \qquad\lim\limits_{t\to\infty}\zeta_-(t)+c^*t=\xi_-.
\end{align}
\end{Theorem}
This will be proved in Section~\ref{sect-convergence.compact.support}. The case of higher dimensions, which involves logarithmic corrections, has been recently studied for the case of the Fisher-KPP reaction term $h(u)=u(1-u)$ in~\cite{Du-Quiros-Zhou-Preprint}.  Such corrections are also expected in higher dimensions for other nonlinearities. This will be considered elsewhere.

In the semilinear case, $m=1$, many similar results exist. When $a=0$, apart from the seminal work of KPP in~\cite{Kolmogorov-Petrovsky-Piscounov-1937}, Stokes introduced the terminology of \emph{pulled} and \emph{pushed} cases ($c^*=\sqrt{2h^\prime (0)}$ and $c^*>\sqrt{2h^\prime (0)}$ respectively) and proved a similar result for the pushed case; see~\cite{Stokes-1976} and also, for more results in the pushed case, the work of Uchiyama~\cite{Uchiyama-1978}. For the pulled case we mention the work of Bramson~\cite{Bramson-1983}. When $a>0$, Kanel~\cite{Kanel-1962} and Fife and McLeod~\cite{Fife-McLeod-1975} proved also similar results. This paper encompass most of these results for the case $m>1$.

\section{A bound for the flux $(u^m)^\prime$}
\label{sect-preliminaries} \setcounter{equation}{0}

In order to identify the asymptotic limits of our solutions as wavefronts, we need to check, following~\cite{Gilding-Kersner-2004}, that these limits satisfy
$$
V\in C(\mathbb{R}), \quad (V^m)^\prime \in L^1_{loc}(\mathbb{R}).
$$
To this aim we need an estimate for the flux $(u^m)^\prime$ of bounded solutions depending only on their size. Such bound will follow from a bound for the derivative of the pressure. In many cases we can get this bound from a semiconvexity result for the pressure, in essence,
$$
\Delta p + \Psi(p) \geq 0,
$$
where $p$ is the pressure related to our solution and $\Psi$ is a function related to our reaction $h$; see~\cite{Perthame-Quiros-Vazquez-2014,Du-Quiros-Zhou-Preprint}. This approach is fine for reactions like $u(1-u)$, but unfortunately we need a different tool due to the many different reactions we are encompassing, since some of them cannot be treated in this way. We will follow the so called Bernstein's method, which requires some effort.

\begin{Lemma}\label{lema_bounded_derivative_pressure}
Let $u$ be a solution of equation~\eqref{eq:main} and suppose that $u_0$ is essentially bounded. The corresponding pressure $p = mu^{m-1}/(m-1)$ satisfies
$$
(p_x)^2 (x,t) \leq \frac{2}{m}\left(\frac{1}{t} + H\right)||p(x,0)||_{\infty},
$$
for some constant $H$ depending only on $\|u_0\|_{\infty}$.
\end{Lemma}

\begin{proof}
First, we define $u_\varepsilon (x,t)$ as the solution of the aproximated problem
$$
u(x,t)=\Delta u^m + h(u), \qquad u(\cdot,0)=u_{0\varepsilon}(\cdot)
$$
where $\{u_{0\varepsilon}\}$ is a family of regular positive functions that converge uniformly in compacts to $u_0$ such that $0<\varepsilon\leq u_{0\varepsilon}\leq M$. We do so to be able to apply the Maximum Principle. It is easy to see that $\varepsilon \leq u_\varepsilon (x,t) \leq M$. Its corresponding pressure satisfies
$$
p_t = (m-1)p\Delta p + p_x^2 + \frac{h(\alpha(p))}{\alpha^\prime(p)}
$$
where $\alpha(p) = ((m-1)p/m)^{1/(m-1)}$ is the density depending on the pressure. We omit the dependence on $\varepsilon$ for the sake of simplicity.

For a fixed $T>0$ we define $S_T=\mathbb{R}\times (0,T]$ and $N=\sup\limits_{x,t \in S_T} p(x,t)$. Now we choose a regular function $\theta$ defined in [0,1] onto $[0,N]$ that is strictly increasing, concave and such that $(\theta^{\prime\prime}/\theta^\prime)^\prime \leq 0$, but keep in mind that a precise description will be specified later. If we define $p=\theta (w)$ our equation becomes
$$
w_t = (m-1)\theta \Delta w + \left( (m-1)\theta \frac{\theta^{\prime\prime}}{\theta^\prime} + \theta^\prime \right) w_x^2 +\frac{h(\alpha(\theta(w)))}{\alpha^\prime \theta^\prime}.
$$

We take now a cutoff function $\zeta \in C(\bar{S_T})\cap C^\infty(S_T)$ such that $0\leq \zeta \leq 1$ and $\zeta \equiv 0$ for $t=0$ or $|x| \geq c >0$ and derive the equation with respect to $x$, multiply it by $w_x\zeta^2$ and consider a point $(x_0, t_0)\in S_T$ where the function $z=w^2_x\zeta^2$ reaches a maximum. In this point we have that $z_t\geq 0,\ z_x=0,\ z_{xx} \leq 0$ and we can discard the case $z\equiv 0$. In other words,
$$
\begin{array}{lcc}
z_t=2w_{xt}w_x \zeta^2 + 2w_x^2\zeta\zeta_t \geq 0 \Rightarrow w_{xt}w_x \zeta^2 \geq -w_x^2\zeta\zeta_t, \\[8pt]
z_x = 2w_{xx}w_x \zeta^2 + 2w_x^2\zeta\zeta_x = 0 \Rightarrow w_{xx}=\frac{-w_x\zeta_x}{\zeta}, \\[8pt]
z_{xx}=2\{ w_{xxx}w_x\zeta^2 + w_{xx}^2\zeta^2 + 4w_{xx}w_x\zeta_x\zeta + w_x^2\zeta_x^2 + \zeta_{xx}w_x^2\zeta\} \leq 0\Rightarrow \\[8pt]
\qquad \Rightarrow w_{xxx}w_x\zeta^2 \leq 2w_x^2\zeta_x^2 - \zeta_{xx}w_x^2\zeta.
\end{array}
$$

Thus
$$
\begin{array}{lcc}
\left( -m\theta^{\prime\prime} - (m-1)\theta\left( \frac{\theta^{\prime\prime}}{\theta^\prime} \right)^\prime \right)\zeta^2 w^4_x \leq -\left((m+1)\theta^\prime +2(m-1)\theta \frac{\theta^{\prime\prime}}{\theta^\prime}\right)\zeta\zeta_x w^3_x \\[8pt]
+ \left( \zeta \zeta_t + 2(m-1)\theta \zeta^2_x - (m-1)\theta\zeta\zeta_{xx} + \zeta^2 \frac{\partial}{\partial w} \left(\frac{h}{\alpha^\prime \theta^\prime} \right) \right) w_x^2.
\end{array}
$$

It is time now to focus on the last derivative that appears in the previous equation. One can check that
$$
\frac{\partial}{\partial w} \left(\frac{h(\alpha)}{\alpha^\prime \theta^\prime} \right) = h^\prime + h \frac{\partial}{\partial w} \left( \frac{1}{\alpha^\prime \theta^\prime} \right) = h^\prime - \frac{h\cdot \alpha^{\prime\prime}(\theta^\prime)^2  + h\cdot \alpha^\prime \theta^{\prime\prime}}{(\alpha^\prime)^2(\theta^\prime)^2}.
$$
All the quantities involved here are bounded, so the only possible problem arises when $\alpha^\prime \to 0$, but bear in mind that $\theta^\prime$ is strictly positive and
$$
\alpha(\theta)=C_1\theta^{1/(m-1)},\quad \alpha^\prime=C_2\alpha^{2-m},\quad \alpha^{\prime\prime} =C_3\alpha^{3-2m}\quad\text{for}\quad C_1,C_2,C_3 \in\mathbb{R},
$$
and thus even when $\alpha^\prime \to 0$ the quantity remains bounded. Therefore, there must exist a positive finite constant $H$, independent of $\varepsilon$ but dependent on $h$ and thus on $\|u_0\|_\infty$, such that
$$
\frac{\partial}{\partial x} \left(\frac{h(\alpha)}{\alpha_\theta \theta_w} \right) \frac{1}{w_x}\leq H.
$$

Let us define
$$
a_1 = \max{|\zeta_t|},\qquad a_2 = \max{|\zeta_x|},\qquad a_3 = \max{|\zeta_{xx}|},
$$
and suppose that there exist some positive constants $b_{1,2,3,4}$ such that
$$
0<Nb_1\leq \theta^\prime \leq Nb_2,\qquad \theta^{\prime\prime} \leq -Nb_3,\qquad | \theta^{\prime\prime}/ \theta^{\prime}|\leq b_4.
$$
Then we have that
$$
\begin{array}{lcc}
(Nmb_3)\zeta^2w_x^4 \leq (N(m+1)b_2 + 2(m-1)Nb_4)a_2\zeta|w_x|^3 \\[8pt]
+(a_1 + 2(m-1)Na_2^2 + (m-1)Na_3 + H)w_x^2,
\end{array}
$$
or, in other words,
$$
\zeta^2w_x^4 \leq c_1(H)w_x^2 + c_2\zeta|w_x|^3.
$$

Since for all $\delta >0$ we have that $c_2\zeta |w_x|^3 \leq \delta \zeta^2w_x^4+  c^2_2w_x^2/4\delta$ (to see this divide by $w_x^2$ and study the resulting parabola), we arrive to
$$
(1-\delta)\zeta^2w_x^4 \leq \left( c_1(H) + \frac{c_2^2}{4\delta} \right)w_x^2
$$
and thus, for all $(x,t)\in S_T$
\begin{equation}\label{bound_z}
z(x,t)\leq \max z \leq \frac{1}{1-\delta} \left(c_1 + \frac{c_2^2}{4\delta}\right).
\end{equation}

This bound depends on the arbitrary functions $\zeta$ and $\theta$, so we fix a point $(x_1,t_1)\in S_T$ and procced to choose them in an appropriate way. We begin with $\zeta$ and choose
$$
\zeta_n(x,t)=\frac{t}{t_1} \psi\left(\frac{x-x_1}{n} \right)
$$
where $\psi$ is a compactly supported regular function satisfying $0\leq \psi\leq 1,\ \psi = 1\ \text{if}\ |x|\leq 1\ \text{and}\ \psi = 0\ \text{if}\ |x|\geq 2$. This way $a_{1n} =1/t_1$ and $a_{2n},a_{3n} \to 0$ as $n \to \infty$, so passing to the limit in~\eqref{bound_z} and substituting $c_1$ and $c_2$ for their values we get
$$
w_x^2(x_1,t_1) = z(x_1,t_1) \leq \frac{1}{Nmb_3t_1} + \frac{H}{Nmb_3}
$$
which means, since $x_1,t_1$ were arbitrary, that
$$
p_x^2(x,t) \leq \frac{Nb_2^2}{mb_3t} + \frac{Nb_2^2H}{mb_3} = \frac{Nb_2^2}{mb_3}\left(\frac{1}{t} + H\right).
$$

On the other hand, we can take $\theta(w) = Nw(a-bw)$ with $a \geq 2b$ and $a \geq b+1$, and minimizing with respect to $a$ and $b$ we get $a=2,\ b=1$. If we recover the subindex $\varepsilon$ we arrive to
$$
(p_{\varepsilon })_x^2(x,t) \leq  \frac{2N}{m}\left(\frac{1}{t} + H\right)
$$
and with this we finish the proof by letting $\varepsilon$ go to 0.
\end{proof}

With this result is easy to prove the following proposition.

\textbf{Proposition 2.1}
Let $u$ be a solution of equation~\eqref{eq:main} with $u_0\in L ^\infty (\mathbb{R})$. Then for every positive time $t >0$ the flux $(u^m)^\prime$ is continuous whenever $u$ remains bounded.

\section{Building sub- and supersolutions}
\label{sect-subsupersolution} \setcounter{equation}{0}

In this section we present some key results that will be needed later following and enhancing the ideas in~\cite{Biro-2002}. They concern the construction of sub- and supersolutions for our problem approaching a TW profile, that we will use all along the rest of the paper. More specifically we study functions of the form
\begin{equation}\label{eq:def_sub_super}
w(x,t)=f(t)V_{c^*}(x-g(t))
\end{equation}
where this $V_{c^*}$ is the unique finite travelling wave solution of our equation. The functions $f$ and $g$ will be asked to solve a system of ODEs in order to guarantee the convergence of $w$ to the travelling wave profile. The reader must see that $f$ represents the height of $w$ while $g$ gives the location of its front.

When we are building the subsolution we want the function $f$ to grow to 1 very slowly so it doesn't exceeds our solution $u$, and we want the function $g$ to increase not to fast for small times so $w$ doesn't surpass $u$ through the boundary while the solution is \lq\lq starting to travel". For big times though, we need $g^\prime (t)\to c^*$ if our intuition is true.

Similarly, when we think about the supersolution we can allow $g$ to grow this time fast at the start while maintaining its behaviour in the limit. This speed also gives a bit more freedom when defining $f$, as we will see, but again it is better to make it go slowly to 1 to ensure that $w$ is above $u$.

Summarizing, our system will be
$$
\left\{
\begin{array}{lcc}
f^\prime (t) =  \varphi(f),\quad f(0)=f_0 \in (1-\delta, 1+\delta), \\[8pt]
g^\prime (t) = c^*f^{m-1} -   k\varphi(f)/f,\quad g(0)=g_0,
\end{array}
   \right.
$$
where $k>0$ is a big enough constant and $\delta$ is such that $h^\prime(u)<0$ for all $u\in(1-\delta,1+\delta)$. This interval must exist due to the continuity of $h^\prime$ in $u=1$.

The function $\varphi:[1-\delta,1+\delta] \rightarrow \mathbb{R}$ is taken such that:
\begin{itemize}
\item $\ \varphi(1)=0,\ \ \varphi^\prime(1)<0,\ \ \varphi(z)>0$ in $[1-\delta,1)$ and $\varphi(z)<0$ in $(1,1+\delta]$.

\item $\ \varphi$ is continuously differentiable (hence locally Lipschitz) in its domain.

\item $\sup\limits_{f\in[1-\delta,1+\delta]} |\varphi^\prime (f)| \leq H(1-\delta)^m$, where $H=\inf\limits_{\theta\in[1-\delta,1+\delta]} |h^\prime (\theta)|$.
\end{itemize}

Based on what we said before, $\varphi$ needs to be a very flat function, almost zero.

\begin{Lemma}\label{lema1}

Consider the following system of ODEs:
\begin{equation}\label{system:sub_super}
\left\{
\begin{array}{lcc}
f^\prime (t) =  \varphi(f),\quad f(0)=f_0 \in (1-\delta, 1+\delta), \\[8pt]
g^\prime (t) = c^*f^{m-1} -   k\varphi(f)/f,\quad g(0)=g_0.
\end{array}
   \right.
\end{equation}
The following properties are satisfied.
\begin{itemize}
\item[\rm (a)] $\lim\limits_{t \rightarrow \infty} f(t) = 1$.

\item[\rm (b)] If $1-\delta < f_0 < 1$ then $f$ is strictly monotone increasing, and hence $g(t) - c^*t$ is strictly monotone decreasing.

\item[\rm (c)]If $1+\delta>f_0 >1$ then $f$ is strictly monotone decreasing, and hence $g(t) - c^*t$ is strictly monotone increasing.

\item[\rm (d)] $\lim\limits_{t \rightarrow \infty} g(t) = \infty$.


\item[\rm (e)] $\lim\limits_{t \rightarrow \infty} g^\prime(t) = c^*$.

\item[\rm (f)] There exists $\xi_0 \in \mathbb{R}$ such that $\lim\limits_{t \rightarrow \infty} (g(t) - c^* t)= \xi_0$.
\end{itemize}
\end{Lemma}

\begin{proof}

Since (a)-(e), existence and uniqueness are obvious, we only have to consider (f). If $f_0=1$, then the result is trivial with $\xi_0=g_0$. If not, then $f^\prime (t) \neq 0$ for all $t>0$ and we can write
$$
\displaystyle\frac{g^\prime(t)- c^*}{f^\prime (t)}= \frac{ c^*(f^{m-1}(t) - 1) - \varphi(f(t)) \cdot (f(t))^{-1}}{\varphi(f(t))}.
$$
Using L'Hopital's rule and the fact that $f \rightarrow 1$ we get
$$
\displaystyle\frac{g^\prime(t) -c^*}{f^\prime (t)} \rightarrow \frac{c^*(m-1)}{(\varphi)^\prime (1)} - k= K
$$
when $t \rightarrow \infty$. Thus  $g^\prime(t) - c^*$ is close to $kf^\prime (t)$ for all large $t$ so there must exist a constant $K^\prime$ such that
\begin{equation}\label{proof:sub_super}
|g^\prime(t) - c^*| \leq K^\prime|f^\prime (t)|\quad \text{for all } t > 0.
\end{equation}

Since $f^\prime$ has a sign, $g^\prime(t) - c^*$ is integrable in $\mathbb{R}^+$, which means that $g(t) - c^*t$ goes to some $\xi_0$ when $t$ grows.
\end{proof}

Let us now check that with the above choices $w$ is indeed a sub- or a supersolution.

\begin{Lemma}\label{lema2}
Let $w$ as in~\eqref{eq:def_sub_super}and $f$ and $g$ as in~\eqref{system:sub_super}. Then:
\begin{itemize}
\item[(i)] If \/ $1-\delta < f_0 < 1$ then $w$ is a subsolution of~\eqref{eq:main}.
\item[(ii)] If \/ $1+\delta>f_0>1$ then $w$ is a supersolution of~\eqref{eq:main}.
\end{itemize}

In both cases there exists $\xi_0 \in \mathbb{R}$ such that
$$
\lim\limits_{t \rightarrow \infty} w(\xi + c^*t, t) = V_{c^*}(\xi - \xi_0)
$$
uniformly with respect to $\xi \in \mathbb{R}$.
\end{Lemma}

\begin{proof}
We start by defining the operator
$$
\mathcal{L}u= u_t - \Delta (u^m) - h(u).
$$

Since $V_{c^*}$ vanishes when $x \geq g(t)$ and the flux for the profile $V_{c^*}$ is continuous, it is enough to study the sign of this operator when $x < g(t)$ to see if $w$ fills the requirements to be a sub/super-solution. Replacing in the equation it is easy to arrive to

\begin{equation}\label{operator}
\mathcal{L}w= \varphi(f)(V_{c^*} + kV_{c^*}^\prime) + f^m h(V_{c^*}) - h(fV_{c^*}).
\end{equation}

\textbf{Remark:} We have made explicit use of the fact that $\Delta(fV_{c^*})^m = f^m\Delta V_{c^*}^m$. This is something that, for example, we can not use in the more general case of the filtration equation, i.e., replacing $u^m$ by a monotone increasing function in the equation.

We want to prove that $\mathcal{L}w\leq 0$, which thanks to~\eqref{operator} is equivalent to \emph{(i),(ii)}.

\noindent\emph{(i)} Now we focus in the case $1-\delta < f_0 < 1$. Let us remark that we can rewrite equation~\eqref{operator} and $\mathcal{L}w\leq 0$ together as
\begin{equation}\label{su.super:ineq}
1+k\frac{V_{c^*}^\prime}{V_{c^*}} + \frac{f^m h(V_{c^*}) - h(fV_{c^*})}{\varphi(f)V_{c^*}} \leq 0.
\end{equation}
We can add and substract $f^m h(fV_{c^*})$ in the last term to get
$$
F_k:=1+k\frac{V_{c^*}^\prime}{V_{c^*}} + \frac{f^m (h(V_{c^*}) - h(fV_{c^*}))}{\varphi(f)V_{c^*}} + \frac{(f^m - 1)h(fV_{c^*})}{\varphi(f)V_{c^*}} \leq 0.
$$
We will check if this is true and study separatedly the cases where $V_{c^*}$ is close to 0 or close to 1.

If we stay where $V_{c^*}\in ((1-\delta)/f_0,1)$, we can apply the Mean Value Theorem in the third term and disregard the second and fourth terms due to their signs to get
$$
F_k\leq 1+\frac{f^m(1-f)h^\prime (\theta)}{\varphi(f)-\varphi(1)}\leq 1+\frac{f^m h^\prime (\theta)}{H(1-\delta)^m} \leq 0.
$$

On the other hand if $V_{c^*}\in [0,(1-\delta)/f_0]$ we can apply again the Mean Value Theorem in both the third and fourth term to see that both those terms are bounded, let us say by $C_1$, and we get
$$
F_k\leq 1+ k\frac{V_{c^*}^\prime}{V_{c^*}} + C_1 \leq 0
$$
if $k>k_1$ for a certain $k_1>0$. Here we have made use of the facts that $V_{c^*}^\prime <0$ in this region and $V_{c^*}^\prime/V_{c^*} \to -\infty$ as $V_{c^*} \to 0$, as mentioned in~\eqref{profile:property}, to ensure that said $k_1$ must exist. This way we finish the part regarding the subsolution when $1-\delta < f_0 < 1$.

\noindent\emph{(ii)} When $1+\delta > f_0 > 1$ we want to see if $\mathcal{L}w\geq 0$ and thus we want again to check~\eqref{su.super:ineq} (notice that $\varphi(f)<0$ this time). Now we add and substract $h(V_{c^*})$ to arrive to
$$
F_k^\prime :=1+k\frac{V_{c^*}^\prime}{V} + \frac{(f^m - 1)h(V_{c^*})}{\varphi(f)V_{c^*}} + \frac{h(V_{c^*})-h(fV_{c^*})}{\varphi(f)V_{c^*}} \leq 0.
$$

First, if $V_{c^*}\in (1-\delta,1)$ then we can, by similar techniques from before, arrive to
$$
F_k^\prime \leq 1+\frac{(1-f)h^\prime (\theta)}{\varphi(f)-\varphi(1)} \leq 1+\frac{h^\prime (\theta)}{H(1-\delta)^m}\leq 0.
$$

If $V_{c^*}\in [0, 1-\delta]$, again by similar techniques, we arrive to
$$
1+ k\frac{V_{c^*}^\prime}{V_{c^*}} + C_2 \leq 0
$$
and this is true if $k>k_2$ for a certain $k_2>0$. We have finished the construction of the supersolutions and for our system of ODEs it is enough to take $k>\max(k_1,k_2)$.

The rest follows from the continuity of $V_{c^*}$ and the existence of a limit for $g(t)-c^*t$.
$$
\lim\limits_{t \rightarrow \infty} w(\xi + c^*t, t) = \lim\limits_{t \rightarrow \infty} f(t)V_{c^*}(\xi + c^*t - g(t)) = \lim\limits_{t \rightarrow \infty}V_{c^*}(\xi + c^*t - g(t)) = V_{c^*}(\xi - \xi_0).
$$
\end{proof}


\section{Propagation of solutions of compact support}
\label{sect-propagation.vs.vanish} \setcounter{equation}{0}

We focus now on initial data that are positive, bounded, piece-wise continuous and compactly supported. We will refer to this as compactly supported initial data.

In this section we will see which conditions on the initial datum and the reaction term $h$ guarantee that the solution will converge to the value $u=1$. This part will advance in parallel to the famous work of Aronson and Weinberger~\cite{Aronson-Weinberger-1978}, in which they develop a similar task for the case $m=1$, the Heat Equation. Many arguments that work out for the Heat Equation are similar (if not simpler, since we don't have to care about the behaviour at infinity) for the PME. The sketch of this part is the following:

First we will study monostable reactions, $a=0$, that are \lq\lq below the Fujita exponent $p_F = m + 2/N$", for whatever this means right now. In this case the solution will propagate no matter the initial datum given it is positive somewhere. This was called, in the semilinear case, the \textit{hair-trigger effect} by Aronson and Weinberger (Chapter 3 in~\cite{Aronson-Weinberger-1978}).

After that we will study at the same time reaction terms with $a=0$ \lq\lq above the Fujita exponent" and reaction terms with $a>0$. In this case the initial datum will be important for it will have to be big enough for us to stablish comparison with a certain subsolution (Chapter 5 in~\cite{Aronson-Weinberger-1978}).

To be able to construct such subsolution we will need to study in advance the travelling wave solutions in the phase-plane of the equation (Chapter 4 in~\cite{Aronson-Weinberger-1978}).

\subsection{Below the Fujita exponent. The hair-trigger effect}

We recall that we are working with $a=0$. Keeping it simple, this term must present a power-like behaviour near the origin, i.e., $h(u)\sim u^p\ \ $ when $\ \ u \sim 0.$

If this exponent $p$ is less or equal than the so called Fujita exponent $p_F$, which in the porous medium equation in dimension $N$ reads
$$
p_F=m + 2/N,
$$
then we will see that our solution converges to 1, regardless of the initial mass of $u_0$. This case is the main focus of this part.

On the other hand, if $p>p_F$ then it is known that there exist solutions with positive initial data that vanish in finite time, as well as others that grow to 1; the phenomena is more complicated. We will later see this case.


The idea goes as follows: we have to show that our solution does not vanish and after that we can prove that it actually converges to 1 uniformly in compact sets. Since we will work in parallel with~\cite{Aronson-Weinberger-1978}, we will only sketch the ideas of the proofs when we consider it necessary.

\begin{Lemma}\label{lema11}
Suppose that our reaction term $h$, apart from the previous hypotheses, satisfies
$$
\displaystyle\liminf\limits_{u \rightarrow 0} \frac{h(u)}{u^{m + 2/N}} > 0.
$$
Then
$$
\displaystyle\limsup\limits_{t \rightarrow \infty} \left\lbrace \sup\limits_{x \in \mathbb{R}} u(x,t) \right\rbrace =1.
$$
\end{Lemma}

The proof of this statement is similar to the one in~\cite{Aronson-Weinberger-1978}. Note that this is not enough to obtain convergence to 1 in compact sets, since our solution can behave like a \lq\lq spike" or converge to 1 in a set that travels to infinity.

The following lemma is based in the work of Kanel'~\cite{Kanel-1964}, and we will present just a sketch of the proof, commenting the principal differences with the one in~\cite{Aronson-Weinberger-1978}.

\begin{Lemma}\label{lema13}

For a fixed $\delta > 0$ define the function
$$
\begin{array}{lcc}
		q_\delta (r) = \delta (1-r^2)^3\quad \text{if}\quad 0 \leq r \leq 1,\\
		\\ q_\delta (r) = 0\quad \text{if}\quad r>1.
\end{array}
$$

Let $v(x,t)$ denotes the solution of the equation
$$
v_t = \Delta (v^m) + \varphi(v),\quad v(x,0) = q_\delta (|x|),
$$
with $0<\delta < \min (b, (3N/k)^{N/2})$. Suppose that $\varphi(v)$ satisfies $a=0$ and
$$
\varphi(v)=k\ v^{m +2/N}\ \  for \ v \in [0,b]
$$
where $k$ is a positive constant and $b \in (0,1)$. Then
$$
\lim\limits_{t \rightarrow \infty} v(0,t)  =1.
$$

\end{Lemma}

\begin{proof}
We divide this sketch of the proof in two parts.

\textit{Part 1} - First, since the initial data is radially simmetric, it is easy through nowadays well known techniques to show that
$$
0 \leq v(x,t) \leq v(0,t) \leq 1.
$$
\textit{Part 2} - Now the hard part is to show that $v(0,t)$ is, after enough time, a monotonic function of time, thus its limit exists. To see this, we set
$$
z = e^{-lt}v_t,\ \ \ l = \max_{w \in [0,1]} (\varphi)^\prime (w).
$$

This is done in order to control the maximal change that $v_t$ can experiment. This would happen with a solution with no diffusion, a flat solution, which will give raise to the ODE
$$
v_t = \varphi(v),
$$
and then $\varphi(v)$ is aproximated via the Mean Value Theorem by $lv$.

If this is so, then $z$ satisfies the equation
$$
z_t = \Delta(mv^{m-1}z) + ((\varphi)^\prime (v) - l)z
$$
with an initial data
$$
z(x,0) = Z(|x|)= \left\{ \begin{array}{lcc}
             m6\delta^m(1-|x|^2)^{3m-3} \{ 6(m-1)|x|^2 + 4|x|^2(1-|x|^2) \\[8pt]
             \qquad - N(1-|x|^2)^2 + \displaystyle \frac{k\delta^{2/N}}{6}(1-|x|^2)^{6/N + 3} \} \quad \text{if }|x|  \leq 1, \\[8pt]
             0  \quad \text{if }|x| >1.
             \end{array}
   \right.
$$

Note that this equation is formally parabolic in the set of positivity of $v$, hence we can use the Maximum Principles.

It is also easy to see that since $\delta < (3N/k)^{N/2}$, there exists an $r_\delta \in (0,1)$ such that $Z(r)$ is an increasing function in $[0, r_\delta]$, $Z(r) <0$ in $[0, r_\delta)$ and $Z(r) \geq 0$ in $[r_\delta, \infty)$.

Thus, the set
$$
S_0 \equiv \{ t >0 : z(0,t) < 0\}
$$
contains the point $t=0$. We shall show that either $S_0= [0,t_1)$ or $S_0=[0,\infty)$ for some $t_1>0$ in order to see the monotonicity of $v(0,t)$.

\textit{Part 2.1} - First we show that $z \equiv 0$ outside a certain ball in every bounded $t$-interval, but this can be easily done by comparing $v$ with a solution of the equation
$$
w_t = \Delta(w^m)+\lambda w
$$
for some finite $\lambda > 0$, thanks to the finite derivative of $\varphi$ in the origin. The supersolution $w$ has finite speed of propagation, and so does $v$ and then $z$.

\textit{Part 2.1} - Second, we have to show that $S_0$ is an interval. This is the most difficult part of the proof, and it is based in both a compactness argument due to the Strong Maximum Principle and an argument involving  continuous paths joining the points $\{t=0, |x| =0\}$ and $\{t=t^\prime > 0, |x| =0\}$ which needs the continuity of $v$ and $z$. This continuity is given in~\cite{diBenedetto-1983,Ziemer-1982}.

The rest of this argument can be found with detail in~\cite{Aronson-Weinberger-1978}, so we will not reproduce it in its whole here since it works in a similar way for us.

\textit{Part 3} - Having that either $S_0= [0,t_1)$ or $S_0=[0,\infty)$, we see that after some time $z$ and $v_t$ have a fixed sign, so $v(0,t)$ is ultimately a monotonic function of time and
$$
\eta^* = \lim\limits_{t \rightarrow \infty} v(0,t)
$$
exists. Now if $\eta^* \in [0,1)$ then for all $\eta \in (\eta^*, 1)$ there exists a $t_n$ such that $v(0,t) < \eta$ for all $t > t_n$, which according to \textit{Part 1} means that
$$
0 \leq v(x,t) \leq v(0,t) < \eta\ \text{for all } t> t_n
$$
which is a contradiction with Lemma~\ref{lema11}, so $\eta^* = 1$ and the lemma is proved.

\end{proof}

At this point we have enough to prove convergence of our solution to 1 in compact sets, the main result of this part.

\begin{Theorem}
Suppose that our reaction term $h$, apart from the previous hypotheses, satisfies
$$
\displaystyle\liminf\limits_{u \rightarrow 0} \frac{h(u)}{u^{m + 2/N}} > 0.
$$

Then if $u$ is a solution of~\eqref{eq:main} with $u \not\equiv 0$ we have that
$$
\liminf\limits_{t \rightarrow \infty} u(x,t) =1
$$
uniformly in bounded subsets of $\mathbb{R}^N$.
\end{Theorem}

The proof is again analogous to the one in~\cite{Aronson-Weinberger-1978}, one only has to be a bit careful about the bounded support of the solution and the comparison stablished in its set of positivity. Notice that this is not the result presented in the introduction, since we lack the information about the speed of propagation $c^*$. The reader will see how we will be able to add this part once we are finished with this section.

Now we would like to study what happens above the Fujita exponent or with reactions with $a>0$. It is well known that even for monostable reactions there are solutions that vanish in this case. It depends on the \lq\lq size" of the initial datum. For example, an initial datum
$$
u_0(x) \leq a
$$
will certainly vanish.

But prior to any results in this direction, we need to study the phase-plane of the travelling waves of the equation.

\subsection{Phase-plane analysis}

We have already stablished in Section 1 the existence of wavefront solutions via the work of Gilding and Kersner, but this is not the point in this section. Our point is to prove the existence of finite travelling waves that do not connect 1 and 0, but $q$ and 0, for a certain value $q>a$.

We do so in order to be able to construct the forementioned three-parameter family of subsolutions. We advance that the solutions we are seeking in this section will be the \lq\lq legs" that connect the flat part of the subsolution and the ground state 0.

We look for solutions $u(x,t) = q(x\cdot \mu -ct),\ \ \xi =x\cdot \mu -ct $ where $\mu$ is an arbitrary unit vector. If we define
$$
p(\xi)=(q^m)^\prime (\xi)
$$
then our equation transforms into
\begin{equation}\label{tweq}
(q^m)^{\prime\prime} (\xi)+ cq^\prime (\xi) + h(q(\xi)) =0
\end{equation}
and we get the system
$$
q^\prime = p/mq^{m-1}, \qquad p^\prime = -cp/mq^{m-1} - h(q).
$$

But, defining $f(q) = mq^{m-1}h(q)$ and $d\xi = mq^{m-1} \ d\tau\ $ we arrive to the less singular system
$$
q^\prime = p,\qquad p^\prime = -cp - f(q).
$$
We must work then in the plane $(q,p)$.

Now we can see that this is the same system that Aronson and Weinberger studied in~\cite{Aronson-Weinberger-1978}. Thus the proofs will work out in a similar way. In fact, a study about the equivalence between these two phase-planes (and hence between their corresponding equations) was made by Engler in~\cite{Engler-1985}.

\textbf{Remark:} There is no need to worry about the change of variable $d\xi = mq^{m-1} \ d\tau\ $. It is just a regularization of the, in some sense, former singular axis $q=0$. The reader may argue that this change of variable is creating an \lq\lq illusory" compact support for trajectories through the mentioned axis, but this is not the case. All the trajectories that we will need will have finite support. Let us show some calculations supporting this claim.

One can integrate in $d\xi = mq^{m-1} \ d\tau$ to obtain
$$
\xi_1 - \xi_0 = \displaystyle\int_{\tau_0}^{\tau_1} mq^{m-1}(\tau)\ d\tau
$$
and see that $\xi$ is going to remain finite whenever $q$ and $\tau$ are finite. This always happen in our trajectories.

Another point of view comes from integrating by separation of variables in the equation $q^\prime = p/mq^{m-1}$ to get
$$
\xi_1-\xi_0 = m \int_{q_0}^{q_1} \frac{q^{m-1}}{p}\ dq.
$$
As we can see, when $q_0\to 0$ the value $p$ goes to a negative value along our trajectories, and thus the integral above is finite, and so is the support of the solution described by the trajectory. It can be proved, in fact, that even the trajectory through $(0,0)$ is finite; an instance of this is, of course, the travelling wave solution $V_{c^*}$ mentioned in the Introduction, which corresponds to a trajectory connecting the singular points $(1,0)$ and $(0,0)$.

Let us start by commenting some properties of this system. For example, if $c=0$ then a trajectory through $(0,0)$ has to satisfy
$$
\frac{p^2}{2} + \int_0^q mu^{m-1}h(u)\ du =0
$$
and this is only possible if this integral is negative for all small values of $q$. Therefore, there is a chance that there are no trajectories through the mentioned point, it depends on $h$.

Also since the linearization of the system in $(0,0)$ presents the eigenvalue $\lambda =0$, we can't use this tool to study this point.

Finally, we comment that the null-cline $\{p^\prime=0\}$ is a curve connecting $(0,0)$ and $(1,0)$ with formula
$$
p(q)=-f(q)/c,
$$
and again it can be checked that $p^\prime(0)=0$.

Now for a given $c \geq 0, \nu>0$ let $p_c(q;\nu)$ be the only trajectory that connects the points $(0,-\nu)$ and $(q_{c,\nu},0)$ with $q_{c,\nu} \in (0,1]$. If $q_{c,\nu} < 1$ then we define $p_c(q;\nu)=0$ for $q\in [q_{c,\nu},1]$. Since $0<\nu<\mu$ implies $p_c(q;\mu) \leq p_c(q;\nu)\leq 0$ for each $q \in [0,1]$ we can define
$$
p_c(q) = \lim\limits_{\nu \rightarrow 0} p_c(q;\nu).
$$

Let $S = \{ (q,p): 0<q<1, p<0 \}$ and
$$
T_c=S\cap \{ (q,p): 0<q<1, p= p_c(q)\}.
$$
Note that perhaps $T_c=\emptyset$.

Define also $q_c \in (0,1]$ as the value such that $p_c(q)<0$ in $(0,q_c)$ and $p_c(q_c)=0$ in the case $q_c \neq 1$. It follows from the monotone convergence theorem that $p_c(q)$ is a solution of (\ref{tweq}) and $T_c$ is a curve through (0,0). If necessary we define $p_c(q) = 0$ for $q\in [q_c,1]$ if $q_c <1$.

Now we state some of the results already stablished in the referred work. Since the proofs are similar, we will only modify the results or say something about the proofs when needed.

\begin{Lemma}\label{lema_pp_1}
For each $c> 0$ the curve $T_c$ is a trajectory of the system in $S$ through (0,0) and it is extremal in the sense that no other trajectory through (0,0) has points in $S$ below $T_c$. Moreover
$$
p_c(q) \geq (1/c) \min\limits_{u \in[0,1]} f(u) - 2cq
$$
in $[0,1]$ and thre exists a $\rho_c \in (0,1]$ such that
$$
p_c(q) \leq \displaystyle\frac{-cq}{2}
$$
in $q \in [0,\rho_c]$.

\end{Lemma}

\begin{Lemma}\label{lema_pp_2}
Suppose that $c^2 > 4\sigma$ where
$$
\sigma \equiv \sup\limits_{u\in [0,1]} f(u)/u < \infty.
$$

Then $q_c =1$ and $p_c(1) < 0$.

\end{Lemma}

In view of these results, the quantity
$$
c^*=\inf\limits_{c>0} \{ q_c=1, p_c(1) <0 \}
$$
is well defined and satisfies $0\leq c^* \leq 4\sigma$.

\begin{Lemma}\label{lema_pp_3}
If $\max\limits_{q\in [0,1]} \displaystyle\int_0^q f(u)\ du >0$ then $c^*>0$.

\end{Lemma}

\textbf{Remark:} Here we see the reason behind the integral condition in~\eqref{eq:cond.Hosono}, but one can ask now what happens when the quantity defined above is negative. We recall the work of Y. Hosono~\cite{Hosono-1986}. In it we see that the sign of $c^*$ actually matches the sign of $s(m,h) = \int_0^1 h(u)u^{m-1}\ du$. When $s(m,h) = 0$ a stationary wave appears, but this case has many similarities to the case $c^* > 0$.

The case $s(m,h) < 0$ is, on the other hand, more interesting. It gives raise to travelling waves with negative speed that are finite but with a front bounded from below for all times, a quite unique feature. This behaviour is used later in~\cite{Hosono-1986} to establish the behaviour of some solutions that vanish but have a support bounded for all times.

\begin{Lemma}\label{lema_pp_4}
If $0 \leq c < d$ and $T_c$ is not empty, then $p_d(q) < p_c(q)$ in $(0,q_c]$. Also, if $d>0$ then $\lim\limits_{c \rightarrow d} p_c(q)=p_d(q)$. Thus, the family $\{T_c\}$ is continuous in the parameter $c$ whenever $c>0$.
\end{Lemma}

\begin{Lemma}
If $c=c^*$ then there exists an unique trajectory connecting the points $(0,0)$ and $(1,0)$ in the system, and thus this $c^*$ and the travelling wave $V_{c^*}$ described by such trajectory are the same as the ones presented in Section~\ref{sect-introduction}.
\end{Lemma}

The following is the final result in this part, but first let us define the quantity
$$
\gamma_c  = \left\{ \begin{array}{lcc}
            0 \text{ if there exists no curve in } S\text{ through } (0,0), \\
            q_c \text{ if the extremal trajectory } T_c\text{ in } S \text{ through }(0,0)\text{ exists.}
             \end{array}
   \right.
$$

\begin{Lemma}\label{lema_pp_5}
If $c \in (0,c^*)$ then $\gamma_c \in [0,1)$ and for every $\eta \in (\gamma_c, 1)$ the only trajectory through $(\eta,0)$ leaves $S$ through the mentioned point and another point in the negative $p$-axis.
\end{Lemma}

\subsection{Over the Fujita exponent}

We are ready now to prove the propagation of solutions whose initial datum is big enough, even when the reaction term is above the Fujita exponent or has $a>0$.

Let us define, for given parameters $c \in (0,c^*),\ \eta \in (\gamma_c, 1)$  and $\rho> (N-1)/c$, the function
$$
v_0(|x|)  = \begin{cases}
             \eta,\quad& |x|\leq \rho, \\
             q(|x| - \rho),\quad& \rho < |x|\leq \rho + b, \\
             0,\quad& |x|\geq \rho + b,
             \end{cases}
$$
where $q(\xi)$ is the trajectory in the previous phase-plane that connects the point $q(0)=(\eta,0)$ with the point $q(b)=(0,-\nu)$, for some $b, \nu >0$.

Let us also define $v$ as the solution for our Cauchy problem with initial datum $v(x,0)=v_0(|x|)$.

\begin{Lemma}\label{lema_propagation}
Under these conditions, we have that
$$
\lim\limits_{t \rightarrow \infty} v(x,t) = 1
$$
uniformly in compact subsets of $\mathbb{R}^N$, and
$$
v(x,t) \geq \eta\quad \text{for}\quad|x| \leq \rho + \left( c - \displaystyle\frac{N-1}{\rho} \right) t\quad \text{and}\quad t \geq 0.
$$

\end{Lemma}

\begin{proof}
Choose an arbitray $c_1\in (0,c - (N-1)/\rho)$ and define
$$
W(x,t) = v_0(|x| - c_1t).
$$

Then we have that

$$
\begin{array}{l}
	W_t(x,t) = -c_1v_0^\prime(|x| - c_1t),\\[8pt]
	\Delta(W^m(x,t)) = \partial_{rr} (v_0^m(r-c_1t)) + \displaystyle\frac{N-1}{|x|} \partial_r(v_0^m(r-c_1t))\
\end{array}
$$
and
$$
h(W(x,t)) = h(v_0(|x| - c_1t)).
$$

Thus we see that
$$
\begin{array}{l}
W_t - \Delta W^m - h(W)  =\\[8pt]
\quad
\begin{cases}
             -h(\eta)\quad&\text{for }|x|\leq \rho + c_1t, \\
             q^\prime(|x| - c_1t - \rho)(c - (N-1)/|x| - c_1) \quad
              &\text{for }\rho + c_1t < |x|\leq \rho + b + c_1t, \\
             0\quad&\text{for }|x|> \rho + b + c_1t.
\end{cases}
\end{array}
$$

From the definition of $q$ we see that the only point where our function is not continuously differentiable is in $|x|=\rho + b +c_1 t$ but this and the fact that $W_t - \Delta W^m - h(W) \leq 0$ are enough to conclude that $W$ is a weak subsolution for our problem, since it can be defined as the supremum of two subsolutions, the function $q$ and 0. Since $W(x,t) \leq v(x,t)\ $ for $x\in \mathbb{R}^N, t\geq 0$ the second assertion of the lemma is proved.

On the other hand, since $v(x,h)\geq W(x,h) \geq W(x,0)=v(x,0)$ for any $h>0$ we have that $v$ is an increasing function of time that is bounded above by the value 1. Thus, $v(x,t)$ converges uniformly on compact sets to a certain $\tau(x) \geq \eta$. At this point it is easy to prove that necessarily $\tau(x) = 1$, and the lemma is proved.

\end{proof}

We can now assert the main result of this section, which is the following.

\begin{Theorem}\label{tma_over_fujita}
Let $u$ be a solution of~\eqref{eq:main} such that
$$
u(x,0) \geq v(x; x_0,\eta,\rho)
$$
for a certain $v$ of the three-parameter family of functions defined in the previous lemma. Then, for any $y_0 \in \mathbb{R}^N$ and $c \in [0,c^*)$
$$
\lim\limits_{t \rightarrow \infty} \min\limits_{|y - y_0| \leq ct} u(y,t) = 1,
$$
anf for any $c>c^*$
$$
\lim\limits_{t \rightarrow \infty} u(y,t) = 0\quad\text{ for }\quad |y - y_0| \geq ct.
$$
\end{Theorem}

\begin{proof}

The first assertion comes from the previous lemma as we see in~\cite{Aronson-Weinberger-1978}. The second one comes from a comparison with a supersolution that we build from the ones presented in Section~\ref{sect-subsupersolution} , $w(x,t)=f(t)V_{c^*}(x-g(t))$. This travelling wave is unidimensional, but we can take, for any unit vector $\nu$ the supersolution defined as $w_\nu(z,t)=w(x\cdot \nu,t)$.

But one has to keep in mind that we can't choose $f_0$ as big as we want, it has to satisfy $1<f_0<1+\delta$ and in principle our solution can surpass this threshold. This is not an issue though, since we can compare with a flat supersolution with initial data $\bar{u}(x,0) \equiv \sup\{u(x,0)\}$.

As we saw in Section~\ref{sect-subsupersolution}, this supersolution $w_\nu$ travels with speed $c^*$ and has a finite front, and thus, if we move with speed greater than $c^*$  in the direction $\nu$ we will only see the value 0, but this $\nu$ was arbitrary, thus we are finished.

\end{proof}
The results presented in this section give convergence on compact sets to the value~1, and also show that the speed of propagation of the solution is $c^*$. But notice that we are lacking the uniform convergence in the whole space.

In other words, we don't know yet how the solution behaves near the front, i.e., near the two free boundaries that appear in our problem. Keep in mind that our intial data may have several disjoint sets of positivity with the corresponding free boundaries, but after a long enough time this sets will merge into one.

We will give uniform convergence in the last section.

\section{Convergence for initial values of class $\mathcal{A}$. Dimension $N= 1$}
\label{sect-tw.class.A} \setcounter{equation}{0}

The aim of this section is to prove convergence towards a travelling wave for solutions with initial data in the class $\mathcal{A}$, Theorem~\ref{thm:convergence.class.A}.

Adapting the work of Bir\'{o}~\cite{Biro-2002} to our equation, we will construct a sub- and a supersolution for this problem (Lemma~\ref{lema2}) which will have certain needed properties (Lemma~\ref{lema1}) and converge to the profile $V_{c^*}$. Then via a special class of solutions, the \textit{eternal solutions} (see next subsection) and a non-degeneracy result (Lemma~\ref{lema:non-degeneracy}) we will see that for certain sequences $t_k \rightarrow \infty$ the solution converges to $V_{c^*}$ both in shape and support. Then we will see that in fact, this is true for any sequence $t \rightarrow \infty$ thanks to a stability result, Lemma~\ref{lema3}.

Notice that we say nothing about the case $\liminf\limits_{x\to -\infty} u_0(x) \leq a$, since there might be points were $u(x,t) > a$ that compensate for the lack of mass at infinity. In the previous section we gave a sufficient condition for the solution to grow to 1 even if it is of compact support, attending only to the concentration of the mass.

In the next lemma we show stability for the solutions to the equation. Essentially, if they start close to $V_{c^*}$, then they remain close to it in some sense.

\begin{Lemma}\label{lema3}
Let $\tilde{u}(\xi, t)$ be a nonnegative continuous solution of the equation~\eqref{eq:main} such that $\tilde{u}(\xi, T) \equiv 0$ for $\xi \geq \overline{\xi}$, $|\overline{\xi}| < \delta$ and $|\tilde{u}(\xi, T) - V_{c^*}(\xi)| < \delta$ for small $\delta > 0$ and fixed $T>0$. Then there exist $\varepsilon(\delta) \to 0$ as $\delta \to 0$, $ \sigma_1(\varepsilon) \geq 0$ and $ \sigma_2(\varepsilon) \geq 0 $ such that for all $t>T$
$$
(1-\varepsilon)V_{c^*}(\xi - \xi_0 + \sigma_1(\varepsilon)) < \tilde{u}(\xi,t) < (1+\varepsilon)V_{c^*}(\xi - \xi_0 - \sigma_2(\varepsilon))
$$
with
$$
  \begin{array}{lcc}
             \displaystyle 0 \leq \sigma_1(\varepsilon) \leq K\varepsilon + \log \frac{1}{1-\varepsilon} \leq K\varepsilon + \frac{\varepsilon}{1-\varepsilon},\\[8pt]
             0 \leq \sigma_2(\varepsilon) \leq K\varepsilon + \log(1+\varepsilon) \leq K\varepsilon + \varepsilon,
             \end{array}
$$
where this $K$ is as in formula~\eqref{proof:sub_super} in the proof of Lemma~\ref{lema1}.

\end{Lemma}

This lemma is a combination of other two that can be found in the work of Bir\'{o} but for his particular reaction term. Since the proof of both is quite similar to the one needed in our case (in fact the proof of what he calls Lemma 2.4 is the same), we refer the interested reader to that paper~\cite{Biro-2002}.

Right now it may seem that we have enough with the previous lemma and the sub- and supersolutions of section 2 for proving the stability of the convergence of our solution, but note that in order to use Lemma~\ref{lema3} we need not only the two functions $\tilde{u}$ and $V_{c^*}$ to be close, but also their supports, which in principle may not happen.

To see this, imagine that our solution converges uniformly to $V_{c^*}$, i.e. to 0 for $\xi \geq \xi^*$ but there is a \textit{\lq\lq thin tail"} in this set that becomes smaller and smaller but with constant (or oscillating) support. In this case we cannot say that the limit of the supports of our solution matches the support of the limit. If this happens we will say that the solution degenerates.

Note that this may only happen in the set $x \geq \xi^*$, since in the complement the uniform convergence of the solution to $V_{c^*}$ impose the convergence of the supports.

Our goal now is to be arbitrarily close to a solution of our problem, both in density and support. In order to fulfill this, we need two things. First, we need to prove that along a certain time sequence our solution does not degenerate. Second we need two auxiliary Cauchy-Dirichlet problems.

\begin{Lemma}\label{lema:non-degeneracy}
For every $r>0$ we have that
$$
\sigma(r) = \limsup\limits_{t\to\infty} u(\zeta(t) - r, t)>0
$$
\end{Lemma}
\begin{proof}
Let us first find a supersolution of our equation. Since $h^\prime (0)<\infty$ there must exist a constant $k$ such that $k u\geq h(u)$. The problem with rection term $k u$ can be transformed via the change of variables
$$
w(x,t) = e^{-kt}v(x,t),\quad \tau = e^{k(m-1)t}/k(m-1)
$$
into the Porous Medium Equation in time $\tau$. It is well known that this equation has a family of explicit travelling waves, and if one writes down their formula and traces back the change of variables it can be seen that the function
$$
v(x,t)=\left[ \frac{m-1}{m}\gamma(t)C_1 \left(\frac{C_1}{k(m-1)}\gamma^{m-1}(t) - x + C_2 \right)_+\right]^{\frac{1}{m-1}},\qquad \gamma(t)=e^{kt}
$$
is a weak supersolution of our equation, with arbitrary $C_1>0$ and $C_2\in\mathbb{R}$. Its free boundary satisfies
$$
\rho(t) = C_2 + C_1\frac{e^{k(m-1)t}}{k(m-1)}
$$
and it is also important to see that $v_x(x_0,t)<0$ if $x_0<\rho(t)$ and that $0\leq v(x,t) \leq v(x,t+h)$ for any positive $h$.

We start now our argument and we will use a contradiction argument. Suppose that there exists a $\delta$ such that $\sigma(\delta)=0$. Then by the eventual monotonicity of $u$ we have that
\begin{equation}\label{eq:max.to.0}
\lim\limits_{t\to\infty} \max\limits_{[\zeta(t)-\delta, \zeta(t)]} u(\cdot,t)=0.
\end{equation}

Set $\delta^*=\min\{\delta/3, c^*/3\}$. By~\eqref{eq:max.to.0} there must exist a pair of values $C_1,C_2$ and a  big enough time $t_0$ such that
\begin{equation}\label{eq:ineq.no.degeneracy}
0\leq u(x,t) < v(\zeta(t_0), t_0)\quad \text{for all } x\in[\zeta(t)-\delta, \zeta(t)],\ t\geq t_0.
\end{equation}
and
$$
\rho(t_0)=\zeta(t_0)+\delta^*,\quad \rho(t_0 + 1)<\zeta(t_0) + 2\delta^*.
$$

We claim that $\zeta(t)\leq\rho(t)$ for $t\in[t_0,t_0+1]$, and if this is so then we will have
$$
\zeta(t_0+1)-\zeta(t_0)<\rho(t_0+1)-\zeta(t_0)<2\delta^*\leq 2c^*/3
$$
for an arbitrarily big $t_0$, and thus the average speed of $\zeta(t)$ can be no larger than $2c^*/3$, which is a contradiction with the fact that $c^*$ is the asymptotic speed of propagation of our solution. Therefore, we only have to prove our claim.

Again, let us assume that it is not true. Then since $\rho(t_0)=\zeta(t_0)+\delta^* > \zeta(t_0)$ and $\zeta$ is continuous there must exist a $t_1\in(t_0,t_0+1)$ such that $\rho(t_0+1) > \zeta(t_1)>\rho(t_1)$, and thus for all $t\in[t_0,t_1]$  we have that
$$
\zeta(t) - \zeta(t_0) + \delta^* < \rho(t_0+1) - \zeta(t_0) + \delta^* < 3\delta^*\leq \delta
$$
and thanks to~\eqref{eq:ineq.no.degeneracy}, for such $t$ we have that
\begin{equation}\label{eq:comparison.no.degeneracy}
u(\zeta(t_0)-\delta^*,t)< v(\zeta(t_0), t_0)< v(\zeta(t_0)-\delta^*, t_0)\leq v(\zeta(t_0)-\delta^*, t).
\end{equation}

Choose $R>\max\{\rho(t_0+1,\zeta(t_0+1))\}$, which means that $R>\zeta(t), \rho(t)$ for $t\in[t_0,t_1]$, and compare $u$ and $v$ in $[\zeta(t_0)-\delta^*, R]\times[t_0,t_1]$. By~\eqref{eq:ineq.no.degeneracy},~\eqref{eq:comparison.no.degeneracy} and the definition of $R$ we have that $u\leq v$ on the parabolic boundary, and thus in the whole region. In particular $u(\rho(t_1),t_1)\leq v(\rho(t_1),t_1)=0$. On the other hand $\zeta(t_1)>\rho(t_1)$, which implies that $u(\rho(t_1),t_1)>0$. This contradiction proves our claim and the lemma.
\end{proof}

This proof tries to be a simplified version of the one found in~\cite{Du-Quiros-Zhou-Preprint} for $h(u)=u(1-~u)$. In this following subsection we introduce the useful concept of \textit{eternal solutions} following the ideas in~\cite{Du-Quiros-Zhou-Preprint}, trying to give a simplified version of them.

\subsection{Eternal solutions}

We are going to study solutions of the following equation:
\begin{equation}
U_t=\Delta(U^m)+c_*U^\prime+h(U)
\end{equation}
defined in $\mathbb{R}^2$, not only in $\mathbb{R}\times \{t\geq 0\}$. From now on this solutions will be denoted with $U$. Such solution will be called a \textit{eternal solution}. They are a useful tool to study the asymptotic behaviour; in some sense they allow us to pass to the limit \lq\lq in two times". First a solution $u$ of problem~\eqref{eq:main} converges to an eternal solution $U$, and then this $U$ converges uniformly to a profile $V_{c^*}$ (we will see that it actually is a profile $V_{c^*}$).

Let us study then a good property present in the eternal solutions. This property is what make them useful for us, and is the following.

\begin{Theorem}\label{tma:eternal}
Let $\xi=x-c^*t$ and $U$ be a nonincreasing (in the variable $\xi$) weak solution to
\begin{equation}
\label{eq:eternal}
U_t=\Delta(U^m)+c_*U^\prime+h(U),\quad (\xi,t)\in\mathbb{R}^2
\end{equation}
such that
$$
V_{c^*}(\xi+C)\leq U(\xi,t)\leq V_{c^*}(\xi-C^*)
$$
for some $C >0$. Then there exists a constant $\xi^*\in[-C,C]$ such that $U(\xi,t)=V_{c^*}(\xi - \xi^*)$.
\end{Theorem}

This is, indeed, a very powerful characterization of eternal solutions. We are saying that every solution $U$ that is nonincreasing and is trapped between two profiles is, actually, a profile. Thus, if $u$ converges to an eternal solution that satisfies the hypotheses of Theorem~\ref{tma:eternal}, it actually converges to a stationary profile.

Let us define
$$
\begin{array}{lcc}
	R_*=\sup\{ R : U(\xi,t)\geq V_{c^*}(\xi-R) \text{ for all } (\xi,t)\in\mathbb{R}^2\},\\[8pt]
	R^*=\inf\{ R : U(\xi,t)\leq V_{c^*}(\xi-R) \text{ for all } (\xi,t)\in\mathbb{R}^2\}.
\end{array}
$$
This pair of values $R_*,R^*$ is well defined and finite due to the hypothesis of our theorem, and in order to prove it is enough to prove that indeed $R_*=R^*$. We will direct our efforts in this direction through the following lemmata. Let us define the free boundary of the function $U$ as $\Psi(t)$.

\begin{Lemma}\label{lema:eternal.converge.above}
There exists a sequence $\{s_n\} \subset \mathbb{R}$ such that
$$
U(\xi, t+s_n)\to V_{c^*}(\xi-R^*),\quad \Psi(t + s_n)\to R^*\quad\text{ as }n\to \infty
$$
uniformly for $(\xi,t)$ in compact subsets of $\mathbb{R}^2$.
\end{Lemma}
\begin{proof}
We define
$$
M(\xi)=\inf\limits_{t\in\mathbb{R}} |V_{c^*}(\xi - R^*) - U(\xi, t)|.
$$
\emph{Step 1.} The first step of the proof is proving that for all $\xi\in\mathbb{R}$, $M(\xi)=0$, which we deduce by contradiction. Suppose that there exists a $\xi_0$ such that $M(\xi_0)=2\varepsilon>0$. We will see that this leads to a contradiction with the definition of $R^*$ by improving the upper bound.

We first improve it in $(-\infty,\xi_0]$. Let us define the auxiliary Cauchy-Dirichlet problem
$$
\left\{\begin{array}{lcc}
		W_t=\Delta(W^m)+c_*W^\prime+h(W)\quad \text{for } \xi\leq \xi_0, t>0 ,\\[8pt]
		W(\xi,0)=1\quad \text{for } \xi\leq \xi_0,\\[8pt]
		W(\xi_0,t)=V_{c^*}(\xi_0-R^*)-\varepsilon\quad \text{for } t>0.
	\end{array}
\right.
$$
It can be proved that this $W$ is a monotone non-increasing function both in time and space and, in fact,
$$
\lim\limits_{t\to\infty} W(\xi,t)=V_{c^*}(\xi_0-(R^* - \delta)),\text{ for  a certain }\delta>0 \text{ fixed}.
$$

On the other hand, by the comparison principle, we have that for each $s\in\mathbb{R}$
$$
U(\xi,s)<W(\xi,0),\quad U(\xi_0,s+t)<W(\xi_0,t) \text{ for }t\geq 0
$$
and thus, by the maximum principle,
$$
U(\xi,s+t)\leq W(\xi,t) \text{ for all }s\in \mathbb{R}, t>0, \xi\leq \xi_0,
$$
which is equivalent to
$$
U(\xi,t)<W(\xi,t-s) \text{ for all }s\in \mathbb{R}, t>s, \xi\leq \xi_0.
$$
Thus, making $s\to - \infty$, we get
$$
U(\xi,t)\leq V_{c^*}(\xi_0-(R^* - \delta)),\text{ for } \xi\leq \xi_0, t\in \mathbb{R}.
$$
An equivalent argument, this time with the Dirichlet problem
$$
\left\{\begin{array}{lcc}
		W_t=\Delta(W^m)+c_*W^\prime+h(W)\quad \text{for } \xi\in [\xi_0, R^*], t>0 ,\\[8pt]
		W(\xi,0)=1\quad \text{for } \xi\in [\xi_0, R^*],\\[8pt]
		W(R^*,t)=0\quad \text{for } t>0,\\[8pt]
		W(\xi_0,t)=V_{c^*}(\xi_0-R^*)-\varepsilon\quad \text{for } t>0
	\end{array}
\right.
$$
shows that
$$
U(\xi,t)\leq V_{c^*}(\xi_0-(R^* - \delta)),\text{ for } \xi\in [\xi_0, R^*], t\in \mathbb{R} \text{ and the same } \delta,
$$
and thus
$$
U(\xi,t)\leq V_{c^*}(\xi_0-(R^* - \delta)),\text{ for } \xi\in \mathbb{R}, t\in \mathbb{R},
$$
contradicting the definition of $R^*$, and hence for all $\xi\in\mathbb{R}$, $M(\xi)=0$, in particular for a certain $\xi_0$ fixed.

Note that we have made use of the fact that $s\in \mathbb{R}$, instead of the more common $s>0$. This means that we are making full usage of the characteristics of eternal solutions to prove that we can actually improve our barrier for all times $t\in \mathbb{R}$.

\noindent\emph{Step 2.} Two possibilities arise now. If $U(\xi_0,t) = V_{c^*}(\xi_0-R^*)$ for a finite time $t$ then the Strong Maximum Principle tells us that indeed $U(\xi,t) \equiv V_{c^*}(\xi-R^*)$ for all $(\xi, t) \in \mathbb{R}^2$, and thus the conclusion of the theorem comes trivially.

Suppose on the contrary that for a certain unbounded sequence $s_n$ we have that $U(\xi_0, s_n)\to V_{c^*}(\xi-R^*)$ and define
$$
U_n(\xi,t)\equiv U(\xi,t+s_n).
$$
The regularity results in~\cite{diBenedetto-1983,Ziemer-1982} and the equiboundedness of $U_n$ allow us to use the Ascoli-Arzel\'a Theorem to conclude that along a certain subsequence $U_n$ converges to $\bar{U}$, a solution of equation~\eqref{eq:eternal} with $\bar{U}(\xi,t)\leq V_{c^*}(\xi-R^*)$ and $\bar{U}(\xi_0,0)= V_{c^*}(\xi_0-R^*)$. Hence $\bar{U}\equiv V_{c^*}$ and the theorem is proved because in this case convergence of solutions impose the convergence of the free boundary $\Psi(t)$.
\end{proof}

Notice that the convergence is in compacts of $\mathbb{R}^2$. We do not get this kind of convergence in the next lemma, but it is not necessary. We need convergence  in compacts of $\mathbb{R}\ni \xi$ along a sequence of times; in particular convergence in compacts of $\mathbb{R}^2$ is a stronger result.

\begin{Lemma}\label{lema:eternal.converge.below}
There exists a sequence $\{\tilde{s}_n\} \subset \mathbb{R}$ such that
$$
U(\xi, \tilde{s}_n)\to V_{c^*}(\xi-R_*),\quad \Psi(\tilde{s}_n)\to R_*\quad\text{ as }n\to \infty
$$
uniformly for $\xi$ in compact subsets of $\mathbb{R}$.
\end{Lemma}
\begin{proof}
We can use again the proof of the previous lemma adapted to a lower barrier. In particular, it is enough to employ this time auxiliary Cauchy-Dirichlet problems with initial data $W(\xi,0)=0$ to conclude, as in the previous lemma, that there exists a sequence $s_j$ such that
$$
U(\xi, t+s_j)\to V_{c^*}(\xi-R_*)\text{ as }j\to \infty.
$$
uniformly for $(\psi,t)$ in compact subsets of $\mathbb{R}^2$.

But note that in this case this is not enough to conclude the convergence of the supports due to de aforementioned possible degeneration of the solutions. We need again a non-degeneracy result similar to Lemma~\ref{lema:non-degeneracy}. This time up to a sequence $\tilde{s}_n$ that we will build ourselves.

Fix a $\delta>0$ and choose an arbitrarily big $k\in \mathbb{N}$ satisfying $k>3(R^*-R_*)/c^*$. We want to see that there exists a $t=t(\delta)\in [0,k]$ and a $j=j(\delta)$ such that
$$
|\Psi(t+s_j)-R_*|<\delta.
$$
For this, we suppose the opposite, that for every $t\in [0,k]$ and $j$ that quantity is bigger or equal than $\delta$, and we will arrive to a contradiction. Once we have this, it is enough to take $\delta_n=1/n$ and $\tilde{s}_n=t(\delta_n) + s_{j(\delta_n)}$ to prove our statement.

We define $\tilde{u}(x,t)=U(x-c^*t,t)$ and $\tilde{\zeta}(t)=\Psi(t)+c^*t$ as its free boundary. Clearly, $\tilde{u}$ is a solution to $\tilde{u}_t=\Delta \tilde{u}^m + h(\tilde{u})$ in $\mathbb{R}^2$.

So fix $\delta>0$, $k>3(R^*-R_*)/c^*$ and $\delta^*<\min(\delta, c^*)/3$. Since $U(\xi, t+\tilde{s}_n)\to V_{c^*}(\xi-R_*)$ uniformly locally in $\mathbb{R}^2$ there must exist a $j$ big enough such that
$$
0<\tilde{u}(\xi + c^*(t + s_j),t + s_j)<\frac{1}{k} \text{ for } t\in [0,k], \xi\in[R_*, \Psi(t + s_j)]\supset [R_*, R_*+\delta].
$$
Now define $v$ as the slow travelling wave supersolution of Lemma~\ref{lema:non-degeneracy} and $\rho(t)$ as its free boundary, and see that because of the previous inequality there must exist for all  $t_0\in[0,k-1]$ two constants $C_1$ and $C_2$ such that
\begin{equation}\label{eq:comparing_eternal}
0\leq \tilde{u}(x,t + s_j)\leq v(\tilde{\zeta}(t_0 +s_j), t_0 +s_j) \text{ for } t\in [0,k], x\geq \tilde{\zeta}(t_0 + s_j)-\delta^*
\end{equation}
and
$$
\rho(t_0 + s_j)=\tilde{\zeta}(t_0+s_j)+\delta^*,\quad \rho(t_0 + 1 + s_j)<\tilde{\zeta}(t_0+s_j) + 2\delta^*.
$$
We will compare $\tilde{u}(x, t+s_j)$ and $v(x, t+s_j)$ in $[\tilde{\zeta}(t_0+s_j)-\delta^*, R]\times [t_0, t_1]$ (remember the definitions of $R$ and $t_1$ from Lemma~\ref{lema:non-degeneracy}).

At the initial time $t=t_0$ there are two possibilities. If $x\in[\tilde{\zeta}(t_0 + s_j)-\delta^*,\tilde{\zeta}({t_0 +s_j})]$, then, since $v$ is nonincreasing in space we have from~\eqref{eq:comparing_eternal} that $\tilde{u}(x, t_0+s_j) \leq v{(x, t_0+s_j)}$. If $x\in[\tilde{\zeta}(t_0 + s_j),R]$ we finish by recalling that here $u\equiv 0$.

In the lateral boundary, when $x=R$ both functions are 0, and when $x=\tilde{\zeta}{(t_0+s_j)-\delta^*}$ we have from the properties of $v$ that
$$
\tilde{u}(\tilde{\zeta}(t_0+s_j)-\delta^*, t+s_j) \leq v(\tilde{\zeta}(t_0+s_j), t_0+s_j) \leq v(\tilde{\zeta}(t_0+s_j) - \delta^*, t+s_j).
$$
So we can compare.

At this point we follow the proof of Lemma~\ref{lema:non-degeneracy}. We have again that $\tilde{\zeta}(t+s_j)\leq\rho{(t + s_j)}$ for $t\in[t_0,t_0+1]$, in particular, this and the fact that $\rho(t_0 + 1 + s_j)<\tilde{\zeta}(t_0+s_j) + 2\delta^*$ provides
$$
\tilde{\zeta}(t_0+1+s_j) - \tilde{\zeta}(t_0+s_j) <2\delta^*<\frac{2c^*}{3}.
$$
But this can be done for every $t_0\in[0,k-1]$. In particular, taking $t_0 = k-1, k-2,...,0$ we see that
$$
  \begin{array}{lcc}
     \tilde{\zeta}(k+s_j) - \tilde{\zeta}(k-1+s_j)<\displaystyle\frac{2c^*}{3},\\[8pt]
     \tilde{\zeta}(k-1+s_j) - \tilde{\zeta}(k-2+s_j)<\displaystyle\frac{2c^*}{3},\\[8pt]
     ...,\\[8pt]
     \tilde{\zeta}(1+s_j) - \tilde{\zeta}(s_j)<\displaystyle\frac{2c^*}{3},
             \end{array}
$$
and a telescoping series tells us that
$$
\tilde{\zeta}(k+s_j) - \tilde{\zeta}(s_j)<\displaystyle\frac{2c^*}{3}k.
$$
Using that $\tilde{\zeta}(t)=\Psi(t)+c^*t$ and $R_*<\Psi(t)<R^*$ we get
$$
k<3(R^*-R_*)/c^*
$$
arriving to a contradiction with $k>3(R^*-R_*)/c^*$.
Thus our $t(\delta)$ and $j(\delta)$ exist and we can construct our sequence $\tilde{s}_n$.
\end{proof}

\begin{Lemma}\label{lema:eternal.R}
$R_*=R^*$.
\end{Lemma}
\begin{proof}
Thanks to Lemma~\ref{lema:eternal.converge.above} and the fact that both $V_{c^*}(\xi-R)$ and $U(\xi,t)$ converge uniformly to 1 when $\xi\to-\infty$ we know that for every $\delta>0$ ther exists a $n(\delta)$ large such that
$$
(1-\delta)V_{c^*}(\xi-R+\delta)\leq U(\xi,s_n)\leq (1+\delta)V_{c^*}(\xi-R-\delta) \text{ for } \xi\in \mathbb{R}.
$$
Thanks to Lemma~\ref{lema3} we have that
$$
(1-\delta)V_{c^*}(\xi-R+\sigma_1(\delta))\leq U(\xi,s_n)\leq (1+\delta)V_{c^*}(\xi-R-\sigma_2(\delta)) \text{ for } \xi\in \mathbb{R}.
$$
with
$$
  \begin{array}{lcc}
             \displaystyle 0 \leq \sigma_1(\varepsilon) \leq K\varepsilon + \log \frac{1}{1-\varepsilon} \leq K\varepsilon + \frac{\varepsilon}{1-\varepsilon},\\[8pt]
             0 \leq \sigma_2(\varepsilon) \leq K\varepsilon + \log(1+\varepsilon) \leq K\varepsilon + \varepsilon,
             \end{array}
$$
where $K$ is independant of $\delta$. It follows that
$$
R^*-\sigma_1(\delta)\leq\Psi(s_n)\leq R^*+\sigma_2(\delta),
$$
and similarly
$$
R_*-\sigma_1(\delta)\leq\Psi(\tilde{s}_m)\leq R_*+\sigma_2(\delta).
$$
Thus for any $t_0>s_n, \tilde{s}_m$ we have that $R^*-\sigma_1(\delta)\leq\Psi(t_0)\leq R_*+\sigma_2(\delta)$, which implies that
$$
0\leq R^*-R_*\leq \sigma_1(\delta)+\sigma_2(\delta)
$$
and we get our result simply by making $\delta$ tend to 0.
\end{proof}

\subsection{Proof of Theorem~\ref{thm:convergence.class.A}}

Having all we need, let us prove our main result for this section. The ideas now are simple enough thanks to the previous work. Basically, what we are going to do is to trap our solution $u$ between sub- and supersolutions and take the limit in time along a non degenerating sequence. This limit will satisfy the conditions of Theorem~\ref{tma:eternal}, and thus it will be a profile $V_{c^*}$. We finish by using the stability result, Lemma~\ref{lema3}.
%

\begin{proof}
We start by setting our sub- and supersolutions below and above $u$,in other words, there must exist $w_1(\xi,t)$ and $w_2(\xi,t)$ coming from Lemma~\ref{lema2} such that
$$
w_1(\xi,t)\leq u(\xi,t)\leq w_2(\xi,t)
$$
and these $w_{1,2}$ converge to different profiles $V_{c^*}$.

On the other hand, thanks to Lemma~\ref{lema3} it is enough to prove that there exists a sequence $t_k$ and a $\xi_0\in\mathbb{R}$ such that
$$
\lim\limits_{k\to\infty} u(\xi,t_k) = V_{c^*}(\xi-\xi_0)
$$
uniformly in compact sets of $\mathbb{R}$ with
$$
\lim\limits_{k\to\infty} \zeta(t_k) = \xi_0.
$$

To do so, if we prove that there exists a time sequence $\{t_k\} \to \infty$ without degeneracy for every $\xi<\zeta(t)$, i.e.,
\begin{equation}\label{eq:no.degeneracy.everywhere}
\liminf\limits_{k\to\infty} u(\zeta(t_k) - r, t_k) >0 \text{ for every } r>0,
\end{equation}
then it is easy to finish, since we can find a subsequence $\{t_{k_j}\}$ such that $u$ converges along it to a eternal solution  without degeneracy, thanks to the regularity of the solution, Alexandrov's reflexion principle and the behaviour in the limit of $w_1$ and $w_2$. The details are left to the reader and can be found in~\cite{Du-Quiros-Zhou-Preprint}. Let us prove then~\eqref{eq:no.degeneracy.everywhere}.

First, we see that thanks to Lemma~\ref{lema:non-degeneracy} and our identification of eternal solutions we have that for each $r_0>0$ there exist a $C_0<0$ and a sequence $t_n$ such that
\begin{equation}\label{eq:fixed.r}
\lim\limits_{n\to\infty} u(\zeta(t_n)-r,t_n) = V_{c^*}(r_0+C_0 - r) \text{ in } C_{loc}(\mathbb{R})\ni r.
\end{equation}

Now we find our sequence $t_k$ via a diagonal argument. Take a sequence $r_k \to 0$ and see that by~\eqref{eq:fixed.r} we have that for each $k>1$ there exist a $C_k<0$ and a sequence $t_n^k$ such that as $n\to \infty$,  $t_n^k \to \infty$ and
$$
u(\zeta(t_n^k)-r,t_n^k) \to V_{c^*}(r_k+C_k - r) \text{ uniformly in } r\in[0,k].
$$
Since $ V_{c^*}(r_k+C_k - r)> V_{c^*}(r_k - r)$ for $r\in [r_k,k]$ there must exist a $n_k$ so large that $t^k_{n_k}>k$ and
$$
u(\zeta(t_{n_k}^k)-r,t_{n_k}^k) > V_{c^*}(r_k - r) \text{ for } r\in [r_k,k].
$$
Define $t_k=t_{n_k}^k$. Then clearly
$$
\liminf\limits_{k\to\infty}  u(\zeta(t_k) - r, t_k) \geq  \lim\limits_{k\to\infty} V_{c^*}(r_k - r) = V_{c^*}(- r)>0.
$$
So the statement and the theorem are proved.
\end{proof}

As a final remark, one can think if the hypothesis $u_0\in\mathcal{A}$ can be weakened. In fact, suppose only that $\liminf_{x\to-\infty} u_0>0$ and let us think how we can reach enough height at $-\infty$.

First, if the reaction term $h$ is under the \textit{hair-trigger effect}, then it is trivial to put a monotone non-increasing in space function below $u_0$ that is going to grow to 1, hence pushing our solution into class $\mathcal{A}$. We shall refer to this as Condition (1), and it depends on $h$.

Second, there is a chance that $u_0$ is not in class $\mathcal{A}$ but it is big enough at $-\infty$ for us to apply the same ideas of Lemma~\ref{lema_propagation} but with a non-symmetric function
$$
v_0(|x|)  = \left\{ \begin{array}{lcc}
             \eta,\quad x\leq \rho, \\
             \\ q(|x| - \rho),\quad \rho < x\leq \rho + b, \\
             \\ 0,\quad x\geq \rho + b,
             \end{array}
   \right.
$$
and push again $u$ into  $\mathcal{A}$. We shall refer to this as Condition (2), and it depends on the mass of $u_0$ near $-\infty$.

\begin{Corollary}\label{coro:lim_inf}
Let $u_0$ be non-negative, bounded, piecewise continuous, $u_0(x) \equiv 0$ for all $x \geq x_0,\ x_0 \in \mathbb{R}$ and satisfying either Condition (1) or Condition (2). Then the same results from Theorem~\ref{thm:convergence.class.A} hold.
\end{Corollary}

\section{Convergence for initial values of compact support. Dimension $N= 1$}
\label{sect-convergence.compact.support} \setcounter{equation}{0}


In this last section we study the uniform convergence of solutions of equation~\eqref{eq:main} with compact support in a similar way as we did in the previous section for initial data in the class $\mathcal{A}$, but let us first see a preliminary result about the speed of convergence to 1 that we will need.

\begin{Lemma}
Let $u$ be a solution of equation~\eqref{eq:main} of compact support such that $u$ converges to 1 uniformly in compact sets. Then there exist $\hat c\in(0,c_*)$, $\delta\in(0,1)$ and $M,T_*>0$
such that
\begin{align}
u(x,t)\leq 1+M e^{-\delta t} & \mbox{ for all $|x|\geq 0$ and $t\geq T_*,$}&\label{u<}\\
u(x,t)\ge 1-Me^{-\delta t} & \mbox{ for all $|x|\in[0,\hat c\,t]$ and $t\ge T_*$.}&\label{u>}
\end{align}
\end{Lemma}

A proof for a harder lemma that this one can be found in~\cite{Du-Quiros-Zhou-Preprint}, but their prove works fine in our case. Let us show now the main proof of this section.


\begin{proof}
At this point, the proof is quite straitforward. We focus only in $\overline{\mathbb{R}_+}$, since in the other side is similar. What we want to do is to trap our solution between a sub- and a supersolution of the form $w_i=f_i(t)V_{c^*}(x-g_i(t))$, the ones that we have already studied. Once this is done, we only have to pass to the moving frame $\xi=x-c^*t$ and repeat the argument showed in the proof of Theorem~\ref{thm:convergence.class.A}.

Thanks to the previous lemma, we focus on times greater than a $T_0$ such that for all $t\ge T_0$
$$
1-Me^{-\delta t}\le u(0,t)\le 1+Me^{-\delta t}
$$
for a fixed $M>0, \delta >0$. To place this $w_i$ below and above the solution in time $T_0$ is easy, but keep in mind that we also have to order the functions in the parabolic frontier $\{x\equiv 0\}\times\{t\geq T_0\}$. In order to do so we only need
$$
f_1(t)V_{c^*}(-g_1(t)) \leq 1-Me^{-\delta t},\qquad 1+Me^{-\delta t} \leq f_2(t)V_{c^*}(-g_2(t)),
$$
but since $V_{c^*}$ is less than 1 ($f_1(t)V_{c^*}(-g_1(t))<f_1(t)$) and strictly monotone decreasing ($f_2(t)V_{c^*}(-g_2(t))\geq f_2(t)V_{c^*}(-g_2(T_0)$), this is the same as asking
$$
f_1(t) \leq 1-Me^{-\delta t},\qquad 1+Me^{-\delta t} \leq C f_2(t),
$$
where $C=V_{c^*}(-g_2(T_0))$ is a positive constant. We will explain how to get the first inequality. The second one is similar.

We can choose $f(T_0) < 1-Me^{-\delta T_0}$, and thus we only have to check that $f^\prime (t) \leq  \delta Me^{-\delta t}$. But remember that $f^\prime=\varphi(f)$. If we choose, for a certain $k>0$,
$$
\varphi(f)= k(1-f),\quad f^\prime(t)=ke^{-kt}
$$
we can arrive to $f^\prime (t) \leq  \delta Me^{-\delta t}$ by chosing the right $k$, and keep in mind that in order to find a proper subsolution we were allowed in Lemma~\ref{lema2} to make $\varphi$ as small as we wanted, so it is enough to pick $\varphi(f)\leq k(1-f)$ for a small enough $k$. One can repeat this argument for the supersolution (this time with a big enough $k$), apply the technique of Theorem~\ref{thm:convergence.class.A} to get convergence in $\Omega_1$ and then repeat in $\Omega_2$ to finish the proof.
\end{proof}


\noindent{\large \textbf{Acknowledgments}}

The author would like to thank Fernando Quir\'os, from Universidad Aut\'onoma de Madrid, for his help and advice in many discussions about the topic.
\




\

\noindent\textbf{Addresses:}

\noindent\textsc{A. G\'arriz: } Departamento de Matem\'{a}ticas, Universidad
Aut\'{o}noma de Madrid, 28049 Madrid, Spain. (e-mail: alejandro.garriz@estudiante.uam.es).


\begin{thebibliography}{99}                                                                                               %







\bibitem{Aronson-Weinberger-1978} Aronson, D.\,G.; Weinberger, H.\,F. \emph{Multidimensional nonlinear diffusion arising in population genetics.} Adv. in Math. 30 (1978), no.\,1, 33--76.



\bibitem{Audrito-Vazquez-2017} Audrito, A.\,V\'{a}zquez, J.\,L. \emph{The Fisher-KPP problem with doubly nonlinear diffusion.} J. Differential Equations 263 (2017), no.\,11, 7647--7708.

\bibitem{Audrito-Vazquez-FDE-2017} Audrito, A.\,V\'{a}zquez, J.\,L. \emph{The Fisher-KPP problem with doubly nonlinear \lq\lq fast'' diffusion.} Nonlinear Anal. 157 (2017), 212--248.


\bibitem{Biro-2002} Bir\'o, Z. \emph{Stability of travelling waves for degenerate reaction-diffusion equations of KPP-type.} Adv. Nonlinear Stud. 2 (2002), no.\,4, 357--371.

\bibitem{Bramson-1983} Bramson, M. \emph{Convergence of solutions of the Kolmogorov equation to travelling waves.} Mem. Amer. Math. Soc. 44 (1983), no.\,285.




\bibitem{diBenedetto-1983}  DiBenedetto, E. \emph{Continuity of weak solutions to a general porous medium equation.} Indiana Univ. Math. J. 32 (1983), no.\,1, 83--118.




\bibitem{Du-Quiros-Zhou-Preprint} Du, Y.; Quir\'{o}s, F.; Zhou, M. \emph{Logarithmic corrections in Fisher-KPP problems for the Porous Medium Equation.} Preprint.

\bibitem{Engler-1985}  Engler, H. \emph{Relations between travelling wave solutions of quasilinear parabolic equations.} Proc. Amer. Math. Soc. 93 (1985), no.\,2, 297--302.

Nonlinear Anal. 43, Ser. A: Theory Methods (2001), no. 8,  943--985.


\bibitem{Fife-McLeod-1975} Fife, P.\,C.; McLeod, J.\,B. \emph{The approach of solutions of nonlinear diffusion equations to travelling wave solutions.} Bull. Amer. Math. Soc. 81 (1975), no.\,6, 1076--1078.

\bibitem{Fisher-1937} Fisher, R.\,A. \emph{The wave of advance of advantageous genes.} Ann. Eugenics 7 (1937), 355--369.



\bibitem{Gilding-Kersner-2004} Gilding, B.\,H.; Kersner, R. \emph{ Travelling waves in nonlinear diffusion-convection reaction}.
Progress in Nonlinear Differential Equations and their Applications, 60. Birkh\"{a}user Verlag, Basel. ISBN: 3-7643-7071-8.

\bibitem{Gurney-Nisbet-1975} Gurney, W.\,S.\,C.; Nisbet, R.\,M.  \emph{The regulation of inhomogeneous populations.}  J. Theoret. Biol.  52 (1975), 441--457.
%

\bibitem{Gurtin-MacCamy-1977} Gurtin, M.\,E.; MacCamy, R.\,C. \emph{On the diffusion of biological populations.}
Math. Biosci. 33 (1977), no.\,1--2, 35--49.



\bibitem{Hosono-1986} Hosono, Y. \emph{Traveling wave solutions for some density dependent diffusion equations.} Japan J. Appl. Math. 3 (1986), no.\,1, 163--196.



\bibitem{Kanel-1962} Kanel', Ja.\,I. \emph{Stabilization of solutions of the Cauchy problem for equations encountered in combustion theory.} (Russian) Mat. Sb. (N.S.) 59 (101) (1962) suppl., 245--288.

\bibitem{Kanel-1964} Kanel', Ja.\,I. \emph{Stabilization of the solutions of the equations of combustion theory with finite initial functions.} (Russian) Mat. Sb. (N.S.) 65 (107) (1964), 398--413.



\bibitem{Kolmogorov-Petrovsky-Piscounov-1937} Kolmogorov, A.;  Petrovsky, I.;  Piscounov, N.
 \emph{\'Etude de l\'equation de la diffusion avec croissance de la quantit\'e de
matire et son application \`a un probl\`eme biologique.} Bull. Univ. \'Etat Moscou
(1937), 1--25.







\bibitem{Newman-Sagan-1981} Newman, W.\,I.; Sagan, C. \emph{Galactic civilizations: population dynamics and interstellar diffusion.} Icarus 46 (1981), 293--327.

\bibitem{dePablo-Vazquez-1991}  de Pablo, A.; V\'azquez, J.\,L. \emph{Travelling waves and finite propagation in a reaction-diffusion equation.} J. Differential Equations 93 (1991), no.\,1, 19--61.

\bibitem{Perthame-Quiros-Vazquez-2014} Perthame, B.; Quir\'os, F.;  V\'azquez, J.\,L. \emph{The Hele-Shaw asymptotics for mechanical models of tumor growth}. Arch. Ration. Mech. Anal.  212 (2014), no.\,1, 93--127.


\bibitem{Sacks-1983} Sacks, P.\,E. \emph{The initial and boundary value problem for a class of degenerate parabolic equations.} Comm. Partial Differential Equations 8 (1983), no.\,7, 693--733.

\bibitem{Samarskii-Galaktionov-Kurdyumov-Mikhailov-Book} Samarskii, A.\,A.; Galaktionov, V.\,A.; Kurdyumov, S.\,P.; Mikhailov, A.\,P. \lq\lq Blow-up in quasilinear parabolic equations''. De Gruyter Expositions in Mathematics, 19. Walter de Gruyter \& Co., Berlin, 1995.  ISBN: 3-11-012754-7.

\bibitem{SanchezGarduno-Maini-1994}  S\'anchez-Gardu\~no, F.; Maini, Ph.\,K. \emph{Existence and uniqueness of a sharp travelling wave in degenerate non-linear diffusion Fisher-KPP equations.} J. Math. Biol. 33 (1994), no.\,2, 163--192.

\bibitem{Stokes-1976} Stokes, A. N. \emph{On two types of moving front in quasilinear diffusion.} Math. Biosci. 31 (1976), no.\,3--4, 307--315

\bibitem{Uchiyama-1978} Uchiyama, K. \emph{The behavior of solutions of some nonlinear diffusion equations for large time.} J. Math. Kyoto Univ. 18 (1978), no.\,3, 453--508.


\bibitem{Zeldovich-1948} Ze1'dovich, J.\,B. \emph{Theory of flame propagation.}   Zhur.
Fiz. Khim. USSR 22 (1948), 27--49.


\bibitem{Ziemer-1982} Ziemer, William P. \emph{Interior and boundary continuity of weak solutions of degenerate parabolic equations.} Trans. Amer. Math. Soc. 271 (1982), no. 2, 733--748.
\



\end{thebibliography}
\end{document}